\documentclass[11pt]{article}
\usepackage{graphicx}
\usepackage{color}
\usepackage{amsmath, amsthm, amssymb}
\usepackage[toc,page]{appendix}
\usepackage{showlabels}
\usepackage{subfigure}
\usepackage{ulem}

\newtheorem{theorem}{Theorem}
\newtheorem{remark}{Remark}

\newtheorem{lemma}{Lemma}
\newtheorem{example}{Example}

\numberwithin{equation}{section}

\begin{document}

\title{A fitted scheme for a Caputo initial-boundary value problem
}


\author{J.L. Gracia\thanks{Department of Applied Mathematics, Torres Quevedo Building, Campus Rio Ebro, University of Zaragoza, 50018 Zaragoza, Spain (jlgracia@unizar.es)}        \and
        E. O'Riordan\thanks{School of Mathematical Sciences, Dublin City University, Glasnevin, Dublin 9, Ireland. (eugene.oriordan@dcu.ie)} \and M. Stynes\thanks{Applied and Computational Mathematics Division, Beijing Computational Science Research Center, Haidian District, Beijing 100193, China. (m.stynes@csrc.ac.cn)}
}




\maketitle

\begin{abstract}
In this paper we consider an initial-boundary value problem with a Caputo time derivative of order $\alpha\in(0,1)$. The solution typically exhibits a weak singularity near the initial time and this causes a reduction in the orders of convergence of standard schemes. To deal with this singularity, the solution is computed with a fitted difference scheme on a graded mesh. The convergence of this scheme is analysed using a discrete maximum principle and carefully chosen barrier functions. Sharp error estimates are proved, which show an enhancement in the convergence rate compared with the standard L1 approximation on uniform meshes, and also indicate an optimal choice for the mesh grading. This optimal mesh grading is less severe than the optimal grading for the standard L1 scheme. Furthermore, the dependence of the error on the final time forms part of our error estimate. Numerical experiments are presented which corroborate our theoretical results.

\end{abstract}

\section{Introduction}

Mathematical models incorporating fractional differential equations arise in many branches of science and enginering (see \cite{EL13,MK00} and the references therein). Hence, it is no surprise that the design and analysis of appropriate numerical methods for solving such fractional differential equations has become a topical  research area within the scientific computing community. In the present paper, we present a new numerical method for solving an initial-boundary value problem with a Caputo  time-fractional derivative of order $\alpha \in (0,1)$.

We shall consider the problem
\begin{subequations} \label{prob}
\begin{align}
&( D_t^\alpha +{\cal L} \bigr)  u(x,t) = f(x,t) \label{proba}
\intertext{for $(x,t)\in Q:=(0,l)\times(0,T]$, with}
&{\cal L}u:=-p \, \frac{\partial^2 u}{\partial x^2} +c(x) u, \label{probL}\\
\intertext{and initial-boundary condtions}
&u(0,t)=u(l,t)=0 \text{ for }t\in(0,T], \label{probb} \\
&u(x,0)=\phi(x)\text{ for } x\in[0,l], \label{probc}
\intertext{where $0<\alpha<1$, $p$ is a positive constant,  $c \in C[0,l]$ with $c\ge 0$, $f\in C(\bar Q)$ where $\bar Q: = [0,l]\times[0,T]$,  and $\phi\in C[0,l]$. Assume the compatibility conditions}
\phi (0)&= \phi (l)=0   \label{phi0l}
\end{align}
\end{subequations}
at the corners $(0,0)$ and $(0,l)$ of $\bar Q$; in Section~\ref{sec:decomposition} more regularity conditions will be imposed on the data of the problem. We consider homogeneous boundary conditions in problem~\eqref{prob} without loss of generality.

The Caputo fractional derivative~\cite{Diet10} of order $\alpha$ with $0<\alpha < 1$  is
\begin{equation}\label{Caputo}
D_t^\alpha \omega(x,t) := \frac1{\Gamma(1-\alpha)}\int_{s=0}^t (t-s)^{-\alpha}\, \frac{\partial \omega(x,s)}{\partial s} \,ds  \quad\text{for } 0\le x \le l, \ 0<t\le T.
\end{equation}

Typical solutions $u(x,t)$ of problem~\eqref{prob} exhibit an initial layer at $t=0$ but nonetheless this singular behaviour has been considered in relatively few papers, e.g.,~\cite{BLY10,JLZ14,JZ16,MAF14,MAFN16,SORG17}. For each fixed $x$, the difference $u(x,t)-\phi(x)$ behaves like a multiple of $t^\alpha$ and consequently the time derivatives of $u$ blow up as $t\to 0$.

In~\cite{SORG17} the convergence in the discrete maximum norm of a standard finite difference scheme is analysed while taking into account this singular behaviour of $u$. The scheme uses the standard L1 approximation for $D^\alpha_t u$ on a graded temporal mesh  and a standard finite difference approximation  for  ${\cal L}u$ on a uniform spatial mesh. It is proved that the method converges with order $O(M^{-2}+N^{-\min\{2-\alpha,r\alpha\}})$ where $M$ and $N$ are the spatial and temporal discretization parameters and $r \ge 1$ is the mesh grading; the value of $r=1$ corresponds to a uniform mesh and the larger $r$ is, the more the grid condenses at $t=0$. This error estimate is sharp and implies that one has to take $r \ge (2-\alpha)/\alpha$ to yield   the optimal order of convergence $O(M^{-2}+N^{-(2-\alpha)})$ method. Thus if $\alpha$ is close to zero, then the mesh will be very fine  near $t=0$ and in practice this could cause some difficulties with rounding error.

In the present paper we construct and analyse a new finite difference scheme that needs less severe grading of the optimal mesh. This scheme is designed to reflect the singular behaviour of the solution near $t=0$: it constructs a difference approximation that is exact for the functions $1$ and $t^\alpha$, i.e., it is \emph{fitted}. We shall prove error estimates for this fitted scheme showing that it is more accurate than the L1 scheme when applied to problems having the initial singularity described earlier.

A related idea was used implicitly in~\cite{JZ16} but the finite element method there is more complicated than ours to implement, and  an $L^\infty(L^2)$-error bound is proved in~\cite{JZ16} while our error bounds are in the stronger $L^\infty(L^\infty)$-norm. See Remark~\ref{rem:JinZhouPODpaper} for more information.

The analysis of convergence of our fitted scheme is based on a decomposition of $u-\phi$ into a leading term that is $O(t^\alpha)$ as $t \to 0$ and a smoother remainder~$v$ (defined in Lemma \ref{lem:decomposition}).  We shall combine truncation error estimates with a novel  stability analysis based on a discrete maximum principle and appropriate barrier functions, obtaining sharp error estimates.

Similarly to~\cite{BHY15}, our discrete maximum principle is proved by means of an alternative definition of the Caputo fractional derivative (see~\cite{BHY15,VainikkoAIP14,VainikkoZAA16}):
\begin{align}
D^\alpha_t \omega(x,t)
    & :=\frac1{\Gamma(1-\alpha)} \left\{ \frac{\omega(x,t)-\omega(x,0)}{t^\alpha} \right. \notag \\
    &\left. \qquad +\alpha \int^{t}_{s=0} (t-s)^{-\alpha-1} [\omega(x,t)-\omega(x,s)]\, ds \right\} \label{Caputo2}
\end{align}
for  $0\le x \le l, \, 0<t\le T$. For all functions considered in the current paper the formulas~\eqref{Caputo} and~\eqref{Caputo2} are equivalent, as can be seen from the conditions for this equivalence that are derived in~\cite{VainikkoAIP14,VainikkoZAA16}.

The alternative formulation~\eqref{Caputo2} is helpful in (i) constructing finite-element type  approximations of $D^\alpha_t$ such as our fitted approximation, (ii) giving a truncation error analysis that is different from the analysis in~\cite{SORG17}, and (iii) facilitating the use of discrete barrier functions in our analysis.

Our error estimates show that in the discrete maximum norm the fitted scheme has convergence  rate $O(M^{-2}+N^{-\min\{2-\alpha,2r\alpha\}})$ (which is an improvement on the rate reported above for~\cite{SORG17}). It follows that the optimal choice of mesh grading for the fitted scheme is $r=\max\{1, (2-\alpha)/(2\alpha) \}$. As well as this theoretical result, it is shown experimentally that the fitted scheme performs better than the standard L1 scheme on (i) uniform meshes (i.e., $r=1$) and (ii) graded meshes with the optimal choice of the mesh grading for each scheme.

The superiority of the fitted scheme for $r=1$ is apparent from the theoretical convergence rates. On optimally graded meshes for both schemes, the fitted scheme has smaller errors because it requires less mesh grading and consequently its mesh is not as severely clustered near $t=0$.

In some applications, large values of the terminal time $T$ are used. Thus it is of interest to explore how error bounds depend on $T$, and our analysis of the fitted scheme will do this. Numerical results will demonstrate the sharpness of our theoretical error bound as a function of~$T$.

The structure of the paper is as follows. In Section~\ref{sec:decomposition} we establish the decomposition of $u$ that will be used to analyse the convergence of the fitted scheme. Estimates for the derivatives of the remainder $v$ near $t=0$ and for large $t$ are given. The fitted scheme is constructed in Section~\ref{sec:discrete}  on a general graded mesh. The convergence of this scheme is analysed in Section~\ref{sec:convergence} and we give our main convergence result where error estimates are proved in the discrete maximum norm. Numerical results in Section~\ref{sec:Numerical} for the L1 and fitted schemes show that the fitted scheme is more accurate than the L1 scheme. In addition, they corroborate the sharpness of our error estimates regarding both the orders of convergence and the dependence of the error on the final time $T$.

\emph{Notation:} We use  the standard Hilbert space $L_2(0,l)$ with norm $\|\cdot\|_2$ and inner product $(\cdot, \cdot)$.   In this paper~$C$ denotes a generic constant that  depends on the data of the boundary value problem~\eqref{prob} but is independent of~$T$ and of any mesh used to solve~\eqref{prob} numerically; note that $C$ can take different values in different places.

\section{Decomposition of the continuous solution} \label{sec:decomposition}

This section gives a decomposition of the solution~$u$ of~\eqref{prob} that will be used to analyse the convergence of the fitted scheme.
The regularity of $u$ is discussed in~\cite[Section 2]{SORG17} and it is shown in~\cite[Theorem 2.1]{SORG17} that
\begin{subequations}\label{Regu}
\begin{align}
\left \vert \frac{\partial^k u}{\partial x^k} (x,t) \right \vert  &\le C \quad\text{for }k=0,1,2,3,4, \label{uxderiv}\\
\left \vert \frac{\partial^{\ell} u}{\partial t^{\ell}} (x,t) \right \vert & \le
C (1+t^{\alpha-\ell}) \quad\text{for } \ell=0, 1,2, \label{utderiv}
\end{align}
\end{subequations}
for all $(x,t)\in [0,l]\times(0,T]$, under some regularity conditions on the data problem.

Set
\[
g(x,t) =  f(x,t)-f(x,0)  + \frac{{\cal L}\bigl({\cal L}\phi(x)-f(x,0)\bigr)}{\Gamma (1+\alpha)}\,t^\alpha \ \text{ for }(x,t)\in Q.
\]

\begin{lemma} \label{lem:decomposition}
Assume that $\phi\in C^4[0,l]$, $c\in C^2[0,l]$ and $f,f_x,f_{xx}\in C(\bar Q)$. Assume also that  $f (0,0)= f (l,0)=0$ and $\phi ''(0)= \phi ''(l)=0$. Then the solution of the problem~\eqref{prob} can be decomposed as
\begin{equation}\label{u-decomposition}
u(x,t) = \phi(x) + \frac{f(x,0)-{\cal L}\phi(x)}{\Gamma (1+\alpha)} t^\alpha + v(x,t)
\end{equation}
where the remainder term $v(x,t)$ is the solution of the following problem:
\begin{align*}
&\bigl (D_t^\alpha +{\cal L} \bigr)  v(x,t) =  g(x,t) \text{ for } (x,t)\in Q, 
\intertext{with initial-boundary conditions}
&v(0,t)= v(l,t)=0 \text{ for }t\in(0,T], \quad
v(x,0)=0\text{ for } x\in[0,l]. \\ 
\end{align*}
\end{lemma}
\begin{proof}
Define
\[
v(x,t) := u(x,t) - \phi(x) + \frac{{\cal L}\phi(x)-f(x,0)}{\Gamma (1+\alpha)}\,t^\alpha \text{ for } (x,t)\in \bar Q.
\]
Then \eqref{u-decomposition} is satisfied,  and $v(0,t)= v(l,t)=0$ for  $t\in(0,T]$ follows from \eqref{probb}, \eqref{phi0l} and our hypotheses on $f$ and $\phi''$. Next, $v(x,0)= 0$ for $x\in [0,l]$ is immediate from~\eqref{probc}. Finally, from~\eqref{proba} and $D^\alpha_t t^\alpha = \Gamma(1+\alpha)$ we get
\[
\bigl (D_t^\alpha +{\cal L} \bigr)  v(x,t) =  f(x,t) - {\cal L}\phi(x) - [f(x,0)-{\cal L}\phi(x)]  + \frac{{\cal L}\bigl({\cal L}\phi(x)-f(x,0)\bigr)}{\Gamma (1+\alpha)}\,t^\alpha,
\]
as desired.
\end{proof}

%
\begin{remark}
The  decomposition~\eqref{u-decomposition} can be motivated as follows. Suppose that
\begin{equation}\label{MotivationDecomposition}
u(x,t)=\phi(x)+z(x)t^\alpha+v(x,t)
\end{equation}
for some function~$v$ whose time derivatives are smoother than those of~$u$; this is a reasonable assumption since we know from~\cite{SORG17} that for each $x$, $u(x,t)-\phi(x)$ behaves like a multiple of~$t^\alpha$. Now the argument of~\cite{StyFCAA} shows that for each fixed $x$,
$$
D^\alpha_t v(x,t)\to 0 \text{ as } t\to 0.
$$
Consequently, substituting~\eqref{MotivationDecomposition} into~\eqref{proba} then letting $t$ go to zero, one gets
$$
{\cal L}\phi(x)+z(x) \Gamma(1+\alpha)=f(x,0),
$$
since $v(x,0)\equiv 0$ and hence ${\cal L}v(x,0)\equiv 0$. Solving this equation for $z(x)$ then substituting into~\eqref{MotivationDecomposition} yields~\eqref{u-decomposition}.
\end{remark}

\begin{remark}
The decomposition~\eqref{u-decomposition} and $\vert u(x,t) \vert \le C$ (from \eqref{uxderiv}) for $(x,t)\in \bar Q$ show that the remainder $v$ satisfies
\begin{align*}
\vert v(x,t) \vert
    = \left \vert u(x,t) - \phi(x) - \frac{f(x,0)-{\cal L}\phi(x)}{\Gamma (1+\alpha)} t^\alpha \right \vert \
     \le C (1+t^\alpha).
\end{align*}
Observe that if $u$ is bounded as $t\to\infty$, then $v$ must blow up like $t^\alpha$ as $t\to\infty$ to ensure that~\eqref{u-decomposition} is valid.  Nevertheless, in Theorem~\ref{thm:BoundRemainderTLarge} we shall see that~$v_t$ and~$v_{tt}$ are bounded as $t\to \infty.$
\end{remark}

The rest of this section is devoted to proving bounds on the remainder~$v$. To derive these we use separation of variables and a regularity framework similar to~\cite{SY11}. Thus  let $\{ (\lambda_i, \psi_i): i=1,2,\dots\}$ be the eigenvalues and eigenfunctions for the Sturm-Liouville two-point boundary value problem
\begin{equation}\label{Sturm}
{\cal L}\psi_i = -p \psi _i'' + c\psi _i = \lambda _i \psi _i \text{ on } (0,l), \quad \psi _i(0)=\psi _i(l)=0.
\end{equation}
The eigenfunctions are normalised by requiring $\|\psi_i \|_2=1$ for all $i$ and  it is well known that $\lambda_i>0$ for all~$i$.

Set
\[
D({\mathcal L}^\gamma) = \left\{ \omega \in L_2(0,l): \sum_{i=1}^\infty \lambda_i^{2\gamma}|(\omega ,\psi_i)|^2 < \infty \right\} \text{ for } \gamma \ge 0,
\]
the domain of the fractional power ${\mathcal L}^\gamma$ of the operator $\mathcal L$ (see~\cite{Henry81,SY11}). Also set
\[
\|\omega \|_{{\mathcal L}^\gamma}= \left(\sum_{i=1}^\infty \lambda_i^{2\gamma}|(\omega,\psi_i)|^2 \right)^{1/2}
\]
and $\omega_i(t):=(\omega(\cdot,t),\psi _i(\cdot))$ for $i=1,2,\dots$

In this framework the bounds~\eqref{Regu} are proved in~\cite[Theorem 2.1]{SORG17} under the hypotheses that
\begin{subequations}\label{assum-regu}
\begin{equation} \label{assum-SIAM}
\phi\in D({\mathcal L}^{5/2}), \,  f(\cdot, t)\in D({\mathcal L}^{5/2}), \,
f_t(\cdot, t) \hbox{ and } f_{tt}(\cdot, t) \in D({\mathcal L}^{1/2})
\end{equation}
  for each $t\in (0,T]$, and
\begin{equation}  \label{assum-SIAMb}
\|f(\cdot,t)\|_{{\mathcal L}^{5/2}} + \|f_t(\cdot,t)\|_{{\mathcal L}^{1/2}} +  \min\{1,t^\rho\} \|f_{tt}(\cdot,t)\|_{{\mathcal L}^{1/2}}\le C
\end{equation}
for all $t\in (0,T]$ and some constant $0<\rho<1$, where $C$ is a constant independent of~$t$. In the present paper we assume that the regularity conditions~\eqref{assum-SIAM} and~\eqref{assum-SIAMb}  are satisfied, and that
\begin{align}
 & f_{ttt}(\cdot, t)  \in D({\mathcal L}^{1/2}),    \label{assum-fttt}\\
&\Vert f_t(\cdot,t) \Vert_{{\cal L}^{1/2}} + \Vert f_{tt} (\cdot,t) \Vert_{{\cal L}^{1/2}} + \min\{1,t^\rho\} \Vert f_{ttt} (\cdot,t) \Vert_{{\cal L}^{1/2}} \le C,  \label{hypofts}
\end{align}
\end{subequations}
for all $t\in(0,T]$,  where $\rho$ and $C$ are constants independent of~$t$ with $0<\rho<1$.

The following calculation appears a few times in our analysis. Suppose that $\omega(\cdot,s)\in D({\mathcal L}^{1/2})$  and $\|\omega(\cdot,s)\|_{{\mathcal L}^{1/2}} \le C$ for $0\le s\le T$. Then
\begin{align}
\sum_{i=1}^\infty |\omega_i(s)|
    & \le \left( \sum_{i=1}^\infty \frac1{\lambda_i} \right)^{1/2} \left( \sum_{i=1}^\infty \lambda_i \omega_i^2(s) \right)^{1/2}
    	\le C \|\omega \|_{{\mathcal L}^{1/2}} \quad\text{for }0\le s\le t, \label{fi0}
\end{align}
from a Cauchy-Schwarz inequality and  $\lambda_i \approx i^2$ \cite[p.335]{CouHil53}.

Theorems~\ref{thm:BoundRemainder} and~\ref{thm:BoundRemainderTLarge}  establish bounds on the derivatives of the remainder~$v$ for $t \leq 1$  and $t >1$, respectively. Their proofs require the generalized Mittag-Leffler function~\cite{Pod99},  which is defined by
\[
E_{\gamma,\delta} (r) := \sum _{k=0}^\infty \frac{r^k}{\Gamma (\gamma k +\delta)} \text{ for } \gamma>0, \,\delta\in\mathbb{R} \text{ and } r\in \mathbb{R}.
\]
One has~\cite[Theorem 1.6]{Pod99}
\begin{equation}\label{MLbound}
 |E_{\gamma, \delta}(-r)|\le \frac{C}{1+r} \le C \text{ for } 0 < \gamma < 2, \,\delta\in\mathbb{R}, \text{ and all } r\ge 0.
\end{equation}

%
\begin{theorem} \label{thm:BoundRemainder}
Assume the hypotheses of Lemma~\ref{lem:decomposition} and~\eqref{assum-regu}. Then there exists a constant $C$ such that for all $(x,t)\in [0,l]\times(0,1]$ one has
\begin{align*}
\left \vert \frac{\partial^{\ell} v}{\partial t^{\ell}} (x,t) \right \vert
    & \le
Ct^{2\alpha-\ell} \quad\text{for } \ell=0,1,2. 
\end{align*}
\end{theorem}
\begin{proof}
We imitate the argument in~\cite[Section 2]{SORG17}. A standard separation of variables technique  yields
$$
v(x,t) = \sum _{i=1}^\infty  \Bigl[\int _{s=0}^t s^{\alpha -1} E_{\alpha ,\alpha }(-\lambda _i s^\alpha) g_i(t-s) \ ds \Bigr] \psi _i (x),
$$
where $g_i(t) =(g(\cdot,t),\psi _i(\cdot))$; cf.~\cite[(2.2)]{SORG17} and recall that $v(x,0)=0$ for $x\in[0,l]$.

To estimate~$v_t$ and~$v_{tt}$ write $v=v^{(1)}+v^{(2)}$, where
\begin{align}
v^{(1)}(x,t)
    & := \sum _{i=1}^\infty  \left[\int _{s=0}^t s^{\alpha -1} E_{\alpha ,\alpha }(-\lambda _i s^\alpha) [f_i(t-s)-f_i(0)] \ ds \right] \psi _i (x), \label{v1}\\
v^{(2)}(x,t)
    & := \sum _{i=1}^\infty   \left[ \left ({\cal L}\bigl({\cal L}\phi(x)-f(x,0)\bigr), \psi_i(x) \right) \right. \notag \\
    & \left. \times \frac1{\Gamma (1+\alpha)} \int _{s=0}^t s^{\alpha -1} E_{\alpha ,\alpha } (-\lambda _i s^\alpha) (t-s)^{\alpha} \, ds \right] \psi _i (x). \label{v2}
\end{align}
The correctness of this decomposition can be seen easily by inspecting Lemma~\ref{lem:decomposition}.

First, consider $v^{(1)}$. Integrating~\eqref{v1} by parts and using~\cite[(4.3.1)]{GKMR14} yields
\begin{equation}\label{v1n}
v^{(1)}(x,t) = \sum _{i=1}^\infty  \left[\int _{s=0}^t s^\alpha E_{\alpha ,1+\alpha }(-\lambda _i s^\alpha) f'_i(t-s) \ ds \right] \psi _i (x).
\end{equation}
Recalling~\eqref{MLbound} and using~\eqref{fi0} and the hypothesis $f_t\in D({\cal L}^{1/2})$ from \eqref{hypofts} gives
\begin{equation}
 \int _{s=0}^t s^\alpha \sum _{i=1}^\infty \vert f'_i(t-s) \vert \, ds
       \le C  \int _{s=0}^t s^\alpha  \, ds \
        \le C t^{\alpha+1}.   \label{salphaf}
\end{equation}
Thus
\begin{equation}
\vert v^{(1)}(x,t) \vert  \le \sum _{i=1}^\infty
	\left[\int _{s=0}^t s^\alpha E_{\alpha ,1+\alpha }(-\lambda _i s^\alpha) \vert f'_i(t-s)\vert  \ ds \right]   \le  C, \label{Bv1}
\end{equation}
using $|\psi_i (x)|\le C$ for all $i$ and $x$ \cite[p.335]{CouHil53}, \eqref{MLbound} and \eqref{salphaf}.

We now consider $v^{(1)}_t$. Differentiate \eqref{v1n} term-by-term with respect to $t$ for $(x,t)\in Q$ and again use  $|\psi_i (x)|\le C$  and~\eqref{MLbound} to get
\begin{align}
    &  \left \vert \sum _{i=1}^\infty  \Bigl[t^\alpha E_{\alpha ,1+\alpha }(-\lambda _i t^\alpha) f'_i(0)+ \int _{s=0}^t s^\alpha E_{\alpha ,1+\alpha }(-\lambda _i s^\alpha) f''_i(t-s) \, ds \Bigr] \psi _i (x)  \right \vert \notag \\
    & \hspace{1cm} \le    Ct^\alpha \sum _{i=1}^\infty  \vert f'_i(0) \vert + C \sum _{i=1}^\infty \int _{s=0}^t  s^\alpha   \vert f''_i(t-s)\vert \, ds
    		\notag \\
    &\hspace{1cm} \le Ct^\alpha, \label{v1t0}
\end{align}
where the $f_i$ terms are bounded as in~\eqref{salphaf}.
Thus for each fixed $t\in [0,1]$, the series in~\eqref{v1t0} is absolutely and uniformly convergent for $(x,t)\in [0,l]\times[0,1]$, so it equals $v^{(1)}_t(x,t)$ on $[0,l]\times[0,1]$ and we obtain
\begin{equation} \label{v1t}
|v^{(1)}_t(x,t)| \le C t^\alpha,  \text{ for } (x,t)\in [0,l]\times[0,1].
\end{equation}

Next, we bound $v^{(1)}_{tt}$. Differentiating~\eqref{v1n} twice term-by-term with respect to $t$  and using \cite[(4.3.1)]{GKMR14}, \eqref{MLbound} and $|\psi_i (x)|\le C$ yields
\begin{align}
    &\hspace{-5mm}  \left \vert \sum _{i=1}^\infty \left[ t^{\alpha-1} E_{\alpha,\alpha} (-\lambda_it^\alpha) f'_i(0)+ t^{\alpha} E_{\alpha ,1+\alpha }(-\lambda _i t^\alpha)  f''_i(0) \right. \right. \notag \\
    & \hspace{1cm} \left. \left.+ \int _{s=0}^t s^{\alpha} E_{\alpha ,1+\alpha }(-\lambda _i s^\alpha)  f'''_i(t-s)   \, ds   \right] \psi _i (x) \right \vert \notag \\
    &\le   Ct^{\alpha-1} \sum _{i=1}^\infty \vert f'_i(0) \vert +  Ct^\alpha  \sum _{i=1}^\infty \vert f''_i(0) \vert + C \sum _{i=1}^\infty \int _{s=0}^t  s^\alpha  \vert f'''_i(t-s) \vert \, ds  \label{v1tt0a}
\end{align}
 for $(x,t)\in [0,l]\times(0,1]$.
The bounds
\begin{equation}  \label{v1tt0b}
\sum _{i=1}^\infty \vert f'_i(0) \vert \le C \text{ and }  \sum _{i=1}^\infty \vert f''_i(0) \vert \le C
\end{equation}
for the first two sums in~\eqref{v1tt0a} are obtained by applying~\eqref{fi0}.
It remains to bound the final sum in~\eqref{v1tt0a}. From~\eqref{hypofts} and~\eqref{fi0} one has
$$
\sum _{i=1}^\infty \vert f'''_i(t-s) \vert \le C (t-s)^{-\rho}.
$$
Hence
\begin{align}
 \int _{s=0}^t s^\alpha \sum _{i=1}^\infty \vert f'''_i(t-s) \vert \, ds
    &  \le C \int _{s=0}^t  s^\alpha (t-s)^{-\rho} \, ds  \le C t^{\alpha+(1-\rho)}  \label{v1tt1}
\end{align}
from a standard identity~\cite[Theorem D.6]{Diet10} for the Euler Beta function.
Now use~\eqref{v1tt0b} and~\eqref{v1tt1} in~\eqref{v1tt0a} to get
    \begin{align}
   &\hspace{-5mm}  \left \vert \sum _{i=1}^\infty \left[ t^{\alpha-1} E_{\alpha,\alpha} (-\lambda_it^\alpha) f'_i(0)+ t^{\alpha} E_{\alpha ,1+\alpha }(-\lambda _i t^\alpha)  f''_i(0) \right. \right. \notag \\
    & \hspace{1cm} \left. \left.+ \int _{s=0}^t s^{\alpha} E_{\alpha ,1+\alpha }(-\lambda _i s^\alpha)  f'''_i(t-s)   \, ds   \right] \psi _i (x) \right \vert \notag \\
    &\le   Ct^{\alpha-1} +  Ct^\alpha  +  C t^{\alpha+(1-\rho)}.  \label{v1tt0}
\end{align}
Thus the series in~\eqref{v1tt0} is absolutely and uniformly convergent for $(x,t)\in [0,l]\times[\varepsilon, t]$ for each $\varepsilon\in (0,t]$, so it equals $v^{(1)}_{tt}$ on $[0,l]\times(0,1]$ and one has
\begin{equation} \label{v1tt}
|v^{(1)}_{tt}(x,t)| \le
C (t^{\alpha-1}+t^\alpha+t^{\alpha+(1-\rho)}) \le C t^{\alpha-1}.
\end{equation}

We now consider $v^{(2)}$. From~\eqref{v2} and~\cite[(1.100)]{Pod99},
\begin{equation}\label{v2t0}
v^{(2)}(x,t) =    \sum _{i=1}^\infty   \left[ t^{2\alpha} \left ({\cal L}\bigl({\cal L}\phi(x)-f(x,0)\bigr), \psi_i\right)  E_{\alpha ,1+2\alpha } (-\lambda _i t^\alpha) \right]\psi _i (x).
\end{equation}
Using $|\psi_i (x)|\le C$ and~\eqref{MLbound} gives
\begin{align}
 \vert v^{(2)}(x,t) \vert
     & \le   C t^{2\alpha} \sum _{i=1}^\infty    \left \vert \left ({\cal L}\bigl({\cal L}\phi(x)-f(x,0)\bigr), \psi_i\right) \right \vert \notag \\
     & = C t^{2\alpha} \sum _{i=1}^\infty \left \vert \left ({\cal L}\phi(x)-f(x,0), {\cal L}\psi_i\right) \right \vert \label{v2t0a}
\end{align}
on  integrating by parts and using  $\phi(0)=\phi(l)=0$, $\phi''(0)=\phi''(l)=0$ and $f(0,0)=f(l,0)=0$.
The series in~\eqref{v2t0a} satisfies
\begin{align*}
\sum _{i=1}^\infty \left \vert \left ({\cal L}\phi(x)-f(x,0), {\cal L}\psi_i\right) \right \vert
    & = \sum _{i=1}^\infty \lambda_i \left \vert \left ({\cal L}\phi(x)-f(x,0), \psi_i\right) \right \vert \\
    & = \sum _{i=1}^\infty \lambda_i \left \vert  \left[  \left (\phi(x), {\cal L}\psi_i\right)- \left (f(x,0), \psi_i\right) \right] \right \vert \\
    & = \sum _{i=1}^\infty \lambda_i \left \vert \left[ \lambda_i \left (\phi(x), \psi_i\right)- \left (f(x,0), \psi_i\right) \right] \right \vert,
\end{align*}
as $(\lambda_i,\psi_i)$ are an eigenpair of~\eqref{Sturm}. Now use $\phi\in D({\cal L}^{5/2})$,   $f(\cdot, t)\in D({\cal L}^{3/2})$ (as $f(\cdot, t)\in D({\cal L}^{5/2})$) and $\|f(\cdot,t)\|_{{\mathcal L}^{3/2}}\le C$ in a calculation resembling~\eqref{salphaf}, to obtain
 \begin{equation} \label{SeriesV2}
\sum _{i=1}^\infty \left \vert \left ({\cal L}\bigl({\cal L}\phi(x)-f(x,0)\bigr), \psi_i\right) \right \vert \le C.
\end{equation}
From~\eqref{v2t0a} and~\eqref{SeriesV2}, we have
\begin{equation} \label{Bv2}
\vert v^{(2)}(x,t) \vert \le C t^{2\alpha}.
\end{equation}

Differentiating~\eqref{v2t0} term-by-term with respect to $t$, from \cite[(4.3.1)]{GKMR14}, \eqref{MLbound} and $|\psi_i (x)|\le C$  one can prove that for $\ell=1,2$ one has
\begin{align}
    & \sum _{i=1}^\infty  \left \vert \frac{\partial^\ell} {\partial t^\ell} \left( t^{2\alpha}  \left[ \left ({\cal L}\bigl({\cal L}\phi(x)-f(x,0)\bigr), \psi_i\right)  E_{\alpha ,1+2\alpha } (-\lambda _i t^\alpha) \right]\psi _i (x)   \right) \right \vert \notag \\
    & \hspace{1cm}  = \sum _{i=1}^\infty  \left \vert  \left( t^{2\alpha-\ell}  \left[ \left ({\cal L}\bigl({\cal L}\phi(x)-f(x,0)\bigr), \psi_i\right)  E_{\alpha ,1+2\alpha-\ell } (-\lambda _i t^\alpha) \right]\psi _i (x)   \right) \right \vert \notag \\
    & \hspace{1cm}  \le C t^{2\alpha-\ell} \sum _{i=1}^\infty  \left \vert    \left ({\cal L}\bigl({\cal L}\phi(x)-f(x,0)\bigr), \psi_i\right)  \right \vert  \notag \\
    & \hspace{1cm} \le C t^{2\alpha-\ell}, \label{v2t-tt0}
\end{align}
 where we have used~\eqref{SeriesV2} in the last inequality. Thus the series in~\eqref{v2t-tt0} is absolutely and uniformly convergent for $(x,t)\in [0,l]\times[\varepsilon, 1]$ for each $\varepsilon\in (0,1]$, so it equals $v^{(2)}_t$ and $v^{(2)}_{tt}$ on $Q$ for $\ell=1,2$ respectively; consequently
\begin{equation} \label{v2t-tt}
\left \vert \frac{\partial^\ell v^{(2)}}{\partial t^\ell} (x,t) \right \vert \le C t^{2\alpha-\ell}
	\text{ for }\ell=1,2 \text{ and } (x,t)\in Q.
\end{equation}

Recall that $v = v^{(1)} + v^{(2)}$.
Combining~\eqref{Bv1}, \eqref{Bv2}, \eqref{v1t}, \eqref{v1tt} and \eqref{v2t-tt}, the theorem is proved.
\end{proof}

A comparison of~\eqref{utderiv} and Theorem~\ref{thm:BoundRemainder} shows that~$v$ is better behaved than~$u$ near $t=0$.

Next, the behaviour of $u$ and $v$  is analysed when $t>1$. Some estimates for~$u$ are proved in Theorem~\ref{thm:BoundUTLarge} and are used  to deduce bounds on~$v$ in Theorem~\ref{thm:BoundRemainderTLarge}.

We shall need a sharper estimate than~\eqref{MLbound} in the particular case when the two parameters of the  Mittag-Leffler function are equal. From~\cite[(1.143)]{Pod99}, one has
\begin{equation}
|E_{\gamma, \gamma}(-r)|  \le \frac{C}{1+r^2} \le  \frac{C}{r^2}  \ \text{ for } 0 < \gamma < 2  \text{ and } r\ge 0. \label{MLboundSpecial}
\end{equation}
\begin{theorem}  \label{thm:BoundUTLarge}
Assume~\eqref{assum-regu} and $T>1$. Then the solution $u$ of~\eqref{prob}  satisfies
\[ 
\left \vert \frac{\partial^\ell u}{\partial t^\ell} (x,t) \right \vert \le C, \ \text{ for }  \ell=1,2 \text{ and } (x,t)\in [0,l]\times[1,T].
\]
\end{theorem}
\begin{proof}
In~\cite[(4.29)]{Luchko12}, \cite[(2.11)]{SY11} or \cite[Section 2]{SORG17}, the solution $u$ to problem~\eqref{prob} is given by
\begin{equation}\label{useries}
u(x,t) = \sum _{i=1}^\infty \Bigl[ (\phi,\psi _i) E_{\alpha,1} (-\lambda _i t^\alpha) +J_i(t)\Bigr] \psi _i (x)
\end{equation}
where
\[
J_i(t) := \int _{s=0}^t s^{\alpha -1} E_{\alpha ,\alpha }(-\lambda _i s^\alpha) f_i(t-s) \ ds
	\quad\text{with}\quad  f_i(t) := (f(\cdot,t),\psi _i(\cdot)).
\]
Differentiating term-by-term with respect to $t$ for $(x,t)\in Q$ yields
\begin{align}\label{utseries}
&\sum _{i=1}^\infty \biggl[ -(\phi,\psi _i) \lambda_i t^{\alpha -1}E_{\alpha,\alpha} (-\lambda _i t^\alpha)
		+ t^{\alpha -1} E_{\alpha ,\alpha }(-\lambda _i t^\alpha) f_i(0)   \Biggr.  \notag \\
	&\hspace{10mm} \Biggl. + \int _{s=0}^t s^{\alpha -1} E_{\alpha ,\alpha }(-\lambda _i s^\alpha)  f_i'(t-s) \, ds \biggr] \psi _i (x),
\end{align}
where a standard formula  \cite[p.154, line 2]{Luchko12} is used to differentiate~$E_{\alpha,1}(\cdot)$. We bound each term as in \cite[Section 2]{SORG17}, obtaining for the first two terms
\begin{align*}
 \sum _{i=1}^\infty \left|(\phi,\psi _i) \lambda_i t^{\alpha -1}E_{\alpha,\alpha} (-\lambda _i t^\alpha)  \psi _i (x)\right|
    & \le Ct^{\alpha-1} \le C, \\
 \sum _{i=1}^\infty  \left| t^{\alpha -1} E_{\alpha ,\alpha }(-\lambda _i t^\alpha) f_i(0) \psi_i(x) \right|
	& \le Ct^{\alpha-1}  \le C.
\end{align*}

Consider now the last term in~\eqref{utseries}. From \eqref{fi0} with $\omega=f_t(\cdot,t)$ and the estimate~\eqref{MLboundSpecial}, we get
\begin{align*}
\int _{s=0}^t s^{\alpha -1} \sum _{i=1}^\infty E_{\alpha ,\alpha }(-\lambda _i s^\alpha)  \left|  f_i'(t-s) \right| \, ds
    & \le  C \int _{s=0}^t \frac{s^{\alpha -1}}{1+s^{2\alpha}} \sum _{i=1}^\infty  \left|f_i'(t-s)\right| \, ds  \\
	& \le  C \int _{s=0}^t \frac{s^{\alpha -1}}{1+s^{2\alpha}} \, ds 	\\
	& \le  C \left (\int _{s=0}^1 s^{\alpha -1} \, ds + \int _{s=1}^t s^{-\alpha -1} \, ds\right) \\
    & \le C.
\end{align*}
Thus,
$$
\left \vert \sum _{i=1}^\infty  \biggl[ \int _{s=0}^t s^{\alpha -1} E_{\alpha ,\alpha }(-\lambda _i s^\alpha)  f_i'(t-s) \, ds \biggr] \psi _i (x) \right \vert \le C,
$$
where we have used $\vert \psi_i \vert \le C$.

From the previous estimates for the three terms in~\eqref{utseries}, we infer that this expression equals $u_t(x,t)$ and
$$
\left \vert \frac{\partial u}{\partial t} (x,t) \right \vert \le C\text{ for }  (x,t)\in [0,l]\times[1,T].
$$

A similar argument  proves the theorem for the case $\ell=2$.
\end{proof}
%

%
\begin{theorem} \label{thm:BoundRemainderTLarge}
Assume the hypothesis of Lemma~\ref{lem:decomposition} and~\eqref{assum-regu}.  Assume that $T>1$. Then
\begin{align*}
\left \vert \frac{\partial^{\ell} v}{\partial t^{\ell}} (x,t) \right \vert & \le
C, \quad\text{for } \ell=1,2 \text{ and } (x,t)\in [0,l]\times[1,T]. 
\end{align*}
\end{theorem}
\begin{proof}
This result follows from the decomposition~\eqref{u-decomposition} and Theorem~\ref{thm:BoundUTLarge}.
\end{proof}

\section{The discretisation} \label{sec:discrete}

We shall approximate the solution of problem~\eqref{prob} using a finite difference scheme on a mesh that is uniform in space and graded in time.

Let $M$ and $N$ be positive integers. Set $h=l/M$ and $x_m:=m h$ for $m=0,1,\ldots,M$. Set $t_n= T(n/N)^r$ for $n=0,1,\ldots,N$ where the constant mesh grading $r\ge 1$ is chosen by the user. If $r=1$, then the mesh is uniform. Set $\tau_n = t_n-t_{n-1}$ for $n=1,2\dots, N$. Then
\begin{equation}\label{tauk+1}
\tau_{n+1} =  T\left(\frac{n+1}{N} \right)^r -  T\left(\frac{n}{N} \right)^r \le CTN^{-r} n^{r-1} = Cn^{-1}t_n
\end{equation}
for $n=0,1,\dots, N-1$.

The mesh is $\{(x_m, t_n): m=0,1,\dots, M, \ n=0,1,\dots, N\}$. The approximate solution computed at the mesh point $(x_m,t_n)$  is denoted by~$u^n_m$.

Our fitted scheme is inspired by the decomposition~\eqref{u-decomposition}.  We approximate $g(x,t)-g(x,s)$ in the Caputo derivative~\eqref{Caputo2} by linear combinations of the functions $1, t^\alpha$ on each temporal subinterval.  That is,  for $0\le s\le T$ consider the basis functions
\begin{subequations} \label{basis}
\begin{align}
\phi_0(s)&=\begin{cases}
\displaystyle{\frac{t^\alpha_1-s^\alpha}{t^\alpha_1}} &\text{for }0\le s \le t_1,\\
0 & \text{otherwise,}
\end{cases} \\
\phi_k(s)&=\begin{cases}
\displaystyle{\frac{s^\alpha-t^\alpha_{k-1}}{t^\alpha_k-t^\alpha_{k-1}}} &\text{for } t_{k-1}\le s \le t_k,\\[3mm]
\displaystyle{\frac{t^\alpha_{k+1}-s^\alpha}{t^\alpha_{k+1}-t^\alpha_k}} &\text{for } t_k\le s \le t_{k+1},\\
0 & \text{otherwise,}
\end{cases}
\quad\text{for }k=1,2,\ldots,N-1,
\\
\phi_N(s)&=\begin{cases}
\displaystyle{\frac{s^\alpha-t^\alpha_{N-1}}{t^\alpha_N-t^\alpha_{N-1}}} &\text{for } t_{N-1}\le s \le t_N,\\
0 & \text{otherwise.}
\end{cases}
\end{align}
\end{subequations}
It is straightforward to check that for $0\le s\le t_n$ and any fixed $x$, the approximation $g(x,s)\approx  \sum_{k=0}^n g(x,t_k) \phi_k(s)$ is exact for  if $g(x,s)$ is constant or $g(x,s)$ is a multiple of $t^s$.

In general, for $0\le s\le t_n$, write
\[
g(x,t_n)-g(x,s) = g(x,t_n)\sum_{k=0}^n \phi_k(s) -g(x,s) \approx \sum_{k=0}^n \phi_k(s) [g(x,t_n) - g(x,t_k)].
\]
Using this approximation in the Caputo fractional derivative formula~\eqref{Caputo2}, $D^\alpha_t u(x_m,t_n)$  is discretised  by
\begin{align}
D^\alpha_N u_m^n
    &:= \frac1{\Gamma(1-\alpha)} \left\{ \frac{u_m^n-u_m^0}{t_n^\alpha} \right. \notag \\
    & \left. \qquad +\alpha \int^{t_n}_{s=0} (t_n-s)^{-\alpha-1}
    	\left[\sum_{k=0}^{n-1} \phi_k(s) (u_m^n-u_m^k) \right]\, ds \right\}  \notag\\
    & = \frac1{\Gamma(1-\alpha)} \left\{\frac{u_m^n-u_m^0}{t^\alpha_n} \right. \notag \\
    & \left. \qquad +\alpha (u_m^n-u_m^0) I^0_n + \alpha \sum^{n-1}_{k=1}  (u_m^n-u_m^k) (I^-_{nk}+I^+_{nk}) \right\},
    	\label{CaputoN1}
\end{align}
where we used the definition~\eqref{basis} of the~$\phi_k$, and
\begin{align*}
I^0_n
    & := \int^{t_1}_{s=t_0} (t_n-s)^{-\alpha-1} \phi_0(s) \, ds \\
    &= -\frac{t^{-\alpha}_n}{\alpha}+\frac1{t_1^\alpha}
    	\int^{t_1}_{s=0} s^{\alpha-1} (t_n-s)^{-\alpha} ds, \\
I^-_{nk}
    & := \int^{t_k}_{s=t_{k-1}} (t_n-s)^{-\alpha-1} \phi_k(s) \, ds \\
    &= \frac{(t_n-t_k)^{-\alpha}}{\alpha}-\frac1{t^\alpha_k-t^\alpha_{k-1}}
    	\int^{t_k}_{s=t_{k-1}}s^{\alpha-1} (t_n-s)^{-\alpha}\, ds, \\
I^+_{nk}
    & := \int^{t_{k+1}}_{s=t_k} (t_n-s)^{-\alpha-1} \phi_k(s) \, ds \\
    &= - \frac{(t_n-t_k)^{-\alpha}}{\alpha}+\frac1{t^\alpha_{k+1}-t^\alpha_k}
    	\int^{t_{k+1}}_{s=t_k}s^{\alpha-1} (t_n-s)^{-\alpha}\, ds.
\end{align*}
Here the integrations by parts are valid since the functions $(t_n-s)^{-\alpha}\phi_k(s)$ are absolutely continuous on $[t_{j-1},t_j]$ for $j=1,2,\ldots,n$ and $k=0,1,\dots, n-1$.

Substituting these three expressions into~\eqref{CaputoN1} and gathering terms, one  obtains the following alternative formulation of~\eqref{CaputoN1}:
\begin{equation}\label{CaputoN2}
D^\alpha_N u_m^n = \frac{\alpha}{\Gamma(1-\alpha)} \sum^{n-1}_{k=0} \frac{u_m^{k+1}-u_m^k}{t^\alpha_{k+1}-t^\alpha_k} \int^{t_{k+1}}_{s=t_k}s^{\alpha-1} (t_n-s)^{-\alpha}\, ds,
\end{equation}
which can be recognised as an approximation of the formula~\eqref{Caputo} for  $D^\alpha_tu(x_m,t_n)$.

The diffusion term in \eqref{probL} is approximated by  a standard second-order discretisation:
$$
u_{xx}(x_m,t_n) \approx \delta^2_x u^n_m := \frac{u_{m+1}^n-2u_m^n+u_{m-1}^n}{h^2}\,.
$$
Thus we approximate \eqref{prob} by the discrete problem
\begin{subequations} \label{scheme}
\begin{align}
L_{M,N} u^n_m &:= D^\alpha_N u^n_m  - p \, \delta^2_x u_m^n  + c(x_m) u^n_m =f(x_m,t_n) \label{schemea}\\
  &\hspace{3cm} \text{ for } 1\le m \le M-1, \ 1 \leq n\leq  N;  \notag \\
u^{n}_0 &=  u^{n}_M=0\ \text{ for } 0< n \le N,  \label{schemeb} \\
u^0_m &=\phi(x_m) \ \text{ for } 0\le m \le M. \label{schemec}
\end{align}
\end{subequations}

\begin{remark}
In practice, the weights of the approximation~\eqref{CaputoN2} are computed from
\begin{align*}
\int^{t_{k+1}}_{s=t_k}s^{\alpha-1} (t_n-s)^{-\alpha}\, ds &= \int^{t_{k+1}/t_n}_{s=t_{k}/t_n} s^{\alpha-1} (1-s)^{-\alpha}\, ds \\
&=B(t_{k+1}/t_n; \alpha,1-\alpha)-B(t_k/t_n; \alpha,1-\alpha),
\end{align*}
where $B(z; \gamma,\delta)$ is the incomplete Beta function defined by
$$
B(z; \gamma,\delta):=\int^z_{s=0} s^{\gamma-1} (1-s)^{\delta-1} \, ds \text{ for } \gamma,\delta>0.
$$

\end{remark}

\begin{remark}\label{rem:JinZhouPODpaper}
In~\cite{JZ16} a method based on the L1 approximation and a proper orthogonal decomposition is used to approximate~\eqref{prob}. It incorporates functions that are constructed to approximate the main singularity in the solution, as our method does. But a  fundamental difference between the two approaches is that from the outset we include explicit basis functions that mimic the main singularity, while~\cite{JZ16} starts with general basis functions then modifies them adaptively to approximate this singularity. Both approaches yield an error term that is $O(t^{2\alpha})$, though we work with $L^\infty(L^\infty)$ (see Theorem~\ref{th:cgce} below) while \cite[Theorem 3.6]{JZ16} uses $L^\infty(L^2)$.
\end{remark}

\section{Convergence analysis of the difference scheme} \label{sec:convergence}

The discretisation of the spatial derivative satisfies
\begin{equation} \label{EstDxx}
 \frac{\partial^2 u}{\partial x^2} (x_m,t_n)= \delta^2_x u(x_m,t_n)+O(h^2).
\end{equation}

Section~\ref{sec:truncCaputo} gives truncation error estimates for the discretisation~\eqref{CaputoN2} of the Caputo temporal derivative. In Section~\ref{sec:staberror} we derive stability and error estimates for our fitted scheme.

\subsection{Truncation error estimates for the temporal fractional derivative}\label{sec:truncCaputo}

Comparing the definition (\ref{Caputo2})  of the Caputo fractional derivative  and the derivation of its discretisation~\eqref{CaputoN1}, one sees that the truncation error can be written as follows:
\begin{align}\label{truncation-error}
&\hspace{-5mm}(D^\alpha_N -D^\alpha_t) u(x_m,t_n)  \notag\\
&= 	\frac{\alpha} {\Gamma(1-\alpha)}  \int^{t_{n}}_{s=0} (t_n-s)^{-\alpha-1} \left[u(x_m,s) - \sum_{k=0}^{n} \phi_k(s) u(x_m,t_k) \right]\, ds \nonumber \\
&= 	\frac{\alpha} {\Gamma(1-\alpha)}  \int^{t_{n}}_{s=0} (t_n-s)^{-\alpha-1} ( u- \bar u)(x_m,s)\, ds,
\end{align}
where $\bar w$ denotes the fitted interpolant of any suitable function~$w$:
\[
\bar w (x_m,s) := \sum_{k=0}^{n} \phi_k(s) w(x_m, t_k) \text{ for }s \in [0,t_{n}].
\]

We saw already that our temporal discretisation is exact for the functions~$1$ and $t^\alpha$ when $x$ is fixed. It then follows from the decomposition~\eqref{u-decomposition} of $u$ that

\begin{subequations}\label{Truncation-error}
\begin{align}
(D^\alpha_N -D^\alpha_t) u(x_m,t_n)
    & = (D^\alpha_N -D^\alpha_t) v(x_m,t_n) \notag \\
    & = \frac{\alpha} {\Gamma(1-\alpha)}  \int^{t_{n}}_{s=0} (t_n-s)^{-\alpha-1} ( v- \bar v)(x_m,s)\, ds   \label{T-e1}.
\end{align}
From~\eqref{Caputo}, \eqref{basis} and~\eqref{CaputoN2} one has the following equivalent representations for the truncation error, which will also be used in our analysis:
\begin{equation} \label{T-e2}
(D^\alpha_N -D^\alpha_t) v(x_m,t_n) = \frac{1} {\Gamma(1-\alpha)}  \int^{t_{n}}_{s=0} (t_n-s)^{-\alpha} \frac{\partial ( v- \bar v)}{\partial s} (x_m,s)\, ds.
\end{equation}
\end{subequations}
Set
\begin{subequations}\label{theta}
\begin{align}
\theta_{n,k} &= \frac{\alpha} {\Gamma(1-\alpha)}  \int^{t_{k+1}}_{s=t_k} (t_n-s)^{-\alpha-1} ( v- \bar v)(x_m,s)\, ds  \label{theta1}
\intertext{for $k=0,1,\dots, n-1$. An integration by parts shows easily that one also has}
\theta_{n,k} &=  \frac{1} {\Gamma(1-\alpha)}  \int^{t_{k+1}}_{s=t_k} (t_n-s)^{-\alpha}
	\frac{\partial ( v- \bar v)}{\partial s} (x_m,s)\, ds \label{theta2}
\end{align}
\end{subequations}
for $k=0,1,\dots, n-2$, and then the equivalence of \eqref{T-e1} and \eqref{T-e2} (or a careful integration by parts of~\eqref{theta1}) implies that~\eqref{theta2} is also valid for $k=n-1$.
Thus we now have the truncation error formula
\begin{equation}\label{trunc}
(D^\alpha_N -D^\alpha_t) u(x_m,t_n) = \sum_{k=0}^{n-1} \theta_{n,k},
\end{equation}
where one can use \eqref{theta1} or \eqref{theta2} to evaluate each $\theta_{n,k}$.


Denote the temporal piecewise linear interpolant of any suitable function~$w$ by
\[
w_I(x_m,s) := w(x_m,t_k) + \frac{w(x_m,t_{k+1})-w(x_m,t_k)}{t_{k+1}-t_k}(s-t_k)\text{ for } s \in I_{k}:=[t_k,t_{k+1}].
\]
For $k \geq 1$ and $w \in C^2[t_1,t_n]$, it is well known that
 \begin{align}
\Vert w-  w_I \Vert _{\infty, I_{k}} \leq C\tau _{k+1} ^2 \Vert w_{tt} \Vert _{\infty, I_{k}} 
\text{ and }
\Bigl\Vert \frac{\partial ( w-  w_I)}{\partial t} \Bigr\Vert _{\infty, I_{k}} \leq C  \tau _{k+1} \Vert w_{tt} \Vert _{\infty, I_{k}}, \label{interp-classical}
\end{align}
where \[
\Vert g \Vert _{\infty, I_{k}} := \max _{x \in [0,l], t\in I_k} \vert g(x,t) \vert.\]

For $s\in I_k$ with $k \geq 0$, the fitted interpolation error can be written as
\begin{equation}\label{fittedlinear}
(v-\bar v)(x_m,s) = (v-v_I)(x_m,s) +\frac{v(x_m, t_{k+1})-v(x_m,t_k)}{t^\alpha _{k+1}-t^\alpha _k}[(s^\alpha) _I -s^\alpha].
\end{equation}
The bounding of this  interpolation error will be split into the two cases of $t_{k+1}\le 1$ and $t_{k+1}>1$.
In the case where $t_{k+1}\le 1$,  Theorem~\ref{thm:BoundRemainder} yields
\[
  \frac{\vert v(x_m,t_{k+1})- v(x_m,t_k) \vert}{t^\alpha_{k+1}-t^\alpha_k}
     =   \frac{ \left \vert \int ^{t_{k+1}}_{s=t_k} \frac{\partial v}{\partial s}(x_m,s) \, ds \right \vert}{t^\alpha_{k+1}-t^\alpha_k}
     \le C\, \frac{ (t^{2\alpha}_{k+1}-t^{2\alpha}_k)}{t^\alpha_{k+1}-t^\alpha_k}
     \leq  C\, t^\alpha_{k+1}.
\]
for $k\ge 0$.
 Then it follows from~\eqref{fittedlinear} that \begin{align*}
\vert (v-\bar v)(x_m,s) \vert
& \le C  \tau^2_{k+1} \Vert v_{tt} \Vert _{\infty, I_{k}} + C t_{k+1}^\alpha \tau^2_{k+1}   \Vert (s^\alpha)'' \Vert _{\infty, I_{k}} \\
& \le C \tau^2_{k+1}  t^{2\alpha-2}_{k} +  C \tau^2_{k+1} t_{k+1}^\alpha t_{k}^{\alpha-2} \\
& \le C \tau^2_{k+1}  t^{2\alpha-2}_{k}.
\end{align*}
In the other case of  $t_{k+1}>1$,  Theorem~\ref{thm:BoundRemainderTLarge}  gives
\[
  \frac{\vert v(x_m,t_{k+1})- v(x_m,t_k) \vert}{t^\alpha_{k+1}-t^\alpha_k} \le  C\frac{t_{k+1}-t_k}{t^\alpha_{k+1}-t^\alpha_k} \leq C  t^{1-\alpha}_{k+1} \text{ for } k\ge 0.     \label{boundFTLarge}
\]
Then \eqref{fittedlinear} yields
\begin{align*}
\vert (v-\bar v)(x_m,s) \vert
& = \vert (v-v_I)(x_m,s)\vert  + \left \vert \frac{v(x_m, t_{k+1})-v(x_m,t_k)}{t^\alpha _{k+1}-t^\alpha _k} \right \vert \, \vert  [(s^\alpha) _I -s^\alpha] \vert \\
& \le C  \tau^2_{k+1} \Vert v_{tt} \Vert _{\infty, I_{k}} + C t_{k+1}^{1-\alpha} \tau^2_{k+1}   \Vert (s^\alpha)'' \Vert _{\infty, I_{k}}  \\
& \le C \tau^2_{k+1}  +   C \tau^2_{k+1} t^{1-\alpha}_{k+1} t_{k+1}^{\alpha-2}  \\
& \leq C \tau^2_{k+1}.
\end{align*}

Combining these two bounds on $v- \bar v$ and using an analogous argument to bound the derivative of the interpolation error, we have  the  interpolation error bounds
\begin{align}
\Vert v- \bar v \Vert _{\infty, I_{k}}
    & \leq  C \tau ^2_{k+1} \max \{1, t_{k}^{2\alpha-2}\}
      \leq  C \tau ^2_{k+1} t_{k}^{2\alpha-2} \max \{1,t_{n}^{2-2\alpha}\} \label{inter-error}
\intertext{and}
\Bigl\Vert \frac{\partial (v- \bar v)}{\partial t} \Bigr\Vert _{\infty, I_{k}}
    &\leq  C\tau _{k+1}  \max \{1, t_{k}^{2\alpha-2}\}.\label{inter-error-FD}
\end{align}

In the next Lemma these interpolation bounds will be used to bound the~$\theta_{n,k}$ in the truncation error \eqref{trunc}.

\begin{lemma}\label{lem:Remainder-CaputoTruncError}
 There exists a constant $C$ such that for all $n \geq 1$
\[
\left| (D^\alpha_N -D^\alpha_t) v(x_m,t_n)\right| \le C \max \{t^\alpha_n, t_{n}^{2-\alpha}\}
	\left[n^{-(2-\alpha)}+  n^{-2r\alpha} \right].
\]
\end{lemma}
\begin{proof}
Fix $n \ge 1$. The argument is based on~\eqref{trunc} and \eqref{theta}.

(a) \emph{Bound on $\theta_{n,0}$:} We use \eqref{theta2}. Observe first that, by Theorem~\ref{thm:BoundRemainder} and the definition (\ref{basis}) of the basis functions, one has
\begin{align*}
\Bigl \vert  \frac{\partial  \bar v}{\partial t} (x_m,t) \Bigr \vert
	&= \Bigl \vert  (v(x_m,t_1) -v(x_m,0)) \phi _0'(t) \Bigr \vert \\
	&= \alpha \frac{ t^{\alpha -1}}{t_1^\alpha} \Bigl \vert \int_{ I_1} \frac{\partial  v}{\partial s} (x_m,s) \ ds  \Bigr \vert \\
	&\leq C t^{\alpha -1}t _1 ^{\alpha }.
\end{align*}
Hence, when $n=1$, one gets
\begin{align}\label{bound1}
 |\theta_{1,0}| &\le \Bigl \vert \int^{t_{1}}_{s=0} (t_1-s)^{-\alpha}  \frac{\partial  \bar v}{\partial s} (x_m,s)  \, ds \Bigr \vert
	+ \Bigl \vert \int^{t_{1}}_{s=0} (t_1-s)^{-\alpha}    \frac{\partial  v}{\partial s} (x_m,s)  \, ds \Bigr \vert \nonumber\\
	&\leq Ct _1 ^{\alpha }  \int^{t_{1}}_{s=0} (t_1-s)^{-\alpha}  s^{\alpha -1}    \, ds
		+ \int^{t_{1}}_{s=0} (t_1-s)^{-\alpha}    s^{2\alpha -1}   \, ds   \nonumber \\
	&\leq  C  t_1^\alpha,
\end{align}
where we used Theorem~\ref{thm:BoundRemainder} and a standard identity~\cite[Theorem D.6]{Diet10} for the Euler Beta function.

By~\eqref{trunc} the proof of the lemma is now complete  in the case $n=1$, so we shall assume that $n>1$ in the rest of the proof.

Then $(t_n-s)^{-\alpha} \le (t_n-t_1)^{-\alpha}$ for $s \in I_1$, so for $n>1$ we have
\begin{align}\label{bound2}
 |\theta_{n,0}|
 	& \le  \int^{t_{1}}_{s=0} (t_n-s)^{-\alpha}\Bigl(  \Bigl \vert  \frac{\partial  \bar v}{\partial s} (x_m,s) \Bigr \vert
		+    \Bigl \vert\frac{\partial  v}{\partial s} (x_m,s) \Bigr \vert \Bigr) \, ds  \notag\\
	&\leq Ct_1^{2\alpha}(t_n-t_1)^{-\alpha} \nonumber \\
	&\leq Ct_1^{2\alpha}t_n^{-\alpha}.
\end{align}

(b) \emph{Bound on $\sum_{k=1}^{\lceil n/2 \rceil -1} \theta_{n,k}$:} (Here $\lceil \cdot \rceil$ is the usual ceiling function.)
Now \eqref{theta1} gives
\begin{align}\label{bound3}
\left| \sum_{k=1}^{\lceil n/2 \rceil -1} \theta_{n,k} \right|
	& \le \sum_{k=1}^{\lceil n/2\rceil -1} \int^{t_{k+1}}_{s=t_k} (t_n-s)^{-\alpha-1}| ( v- \bar v)(x_m,s)| \, ds 	\nonumber\\
	& \leq \sum_{k=1}^{\lceil n/2\rceil -1} \Vert v- \bar v \Vert _{\infty, I_{k}}  \int^{t_{k+1}}_{s=t_k}
		(t_n-t_{\lceil n/2 \rceil})^{-\alpha-1}  \, ds \nonumber \\
	& \leq C \sum_{k=1}^{\lceil n/2\rceil -1}\Vert v- \bar v \Vert _{\infty, I_{k}}    t_n^{-\alpha-1} \tau _{k+1} \nonumber \\
	& \leq  C\max \{1, t_{n}^{2-2\alpha}\}  \sum_{k=1}^{\lceil n/2\rceil -1} \Bigl(\frac{\tau _{k+1}}{t_k}\Bigr)^{3 }
		\Bigl(\frac{t _{k}}{t_n}\Bigr)^{1+2\alpha }t_n^\alpha,  \nonumber\\
	&\leq   C \max \{t^\alpha _n, t_{n}^{2-\alpha}\} \sum_{k=1}^{\lceil n/2\rceil -1} k^{-3}
		(k^rn^{-r})^{(2\alpha +1)} t_k, \nonumber\\
	&=  C \max \{t^\alpha _n, t_{n}^{2-\alpha}\} n^{-r(2\alpha +1)}\sum_{k=1}^{\lceil n/2\rceil -1}  k^{r(2\alpha +1) -3},
\end{align}
where we used the bounds (\ref{inter-error}) and \eqref{tauk+1}.
 As in ~\cite[(5.9)]{SORG17}, we can use the inequality
\begin{align}\label{bound4}
n^{-r(2\alpha+1)}\sum_{k=1}^{\lceil n/2\rceil -1}  k^{r(2\alpha +1)-3}
     & \le C
	\begin{cases}
		  n^{-r(2\alpha+1)} &\text{if } r(2\alpha+1)<2, \\
		  n^{-2}\ln n &\text{if } r(2\alpha+1)=2, \\
		  n^{-2} &\text{if } r(2\alpha+1)>2,	
	\end{cases}
\end{align}
to bound the sum in (\ref{bound3}).

(c)  \emph{Bound on $\sum_{k=\lceil n/2 \rceil}^{n-2} \theta_{n,k}$:}  Invoke \eqref{theta1} and \eqref{inter-error},
and observe that  $\tau_k \leq \tau_n$ and $t_n \ge t_k \ge  2^{-r}t_n$ for  $\lceil n/2 \rceil \le k \le n$, to get
\begin{align}\label{bound5}
\left| \sum_{k=\lceil n/2 \rceil}^{n-2} \theta_{n,k} \right|
	&\le C \left \vert\int^{t_{n-1}}_{s=t_{\lceil n/2 \rceil}} (t_n-s)^{-\alpha-1} (v- \bar v)(x_m,s) \, ds \right \vert \nonumber \\
	&\leq C \max \{1, t_{n}^{2-2\alpha}\} \tau _n^{2} t_n^{2\alpha -2} \int^{t_{n-1}}_{s=t_{\lceil n/2 \rceil}}
		(t_n-s)^{-\alpha-1} \, ds   \nonumber \\
	&\leq C\max \{1, t_{n}^{2-2\alpha}\} \tau _n^{2} t_n^{2\alpha -2}  (t_n-t _{n-1}) ^{-\alpha} \nonumber  \\
	& = C \max \{t^\alpha_n, t_{n}^{2-\alpha}\} \Bigl(\frac{\tau _n}{t_n}\Bigr)^{(2-\alpha) } \nonumber \\
	&\leq C \max \{t^\alpha_n, t_{n}^{2-\alpha}\} n^{-(2-\alpha) },
\end{align}
by \eqref{tauk+1}.

(d)  \emph{Bound on $\theta_{n,n-1}$:} By \eqref{theta2} and \eqref{inter-error-FD} we have
\begin{align}\label{bound6}
|\theta_{n,n-1}| &\leq \left\Vert \frac{\partial (v- \bar v)}{\partial t} \right\Vert _{\infty, I_{n-1}}   \int^{t_{n}}_{s=t_{n-1}} (t_n-s)^{-\alpha} \ ds \nonumber \\
 & \leq C\max \{1, t_{n}^{2-2\alpha}\}  \tau _n  t _{n}^{2\alpha -2} (t_n-t _{n-1})^{1-\alpha}\nonumber  \\
& \leq
 C \max \{t^\alpha_n, t_{n}^{2-\alpha}\} n^{-(2-\alpha) },
\end{align}
similarly to \eqref{bound5}.

Collecting the bounds (\ref{bound1})-(\ref{bound6}) yields (for $n >1$),
\begin{align*}
&\hspace{-3mm}\left\vert ( D^\alpha_N -D^\alpha_t) v(x_m,t_n)\right\vert \\
&\leq Ct_1^{2\alpha} t_n ^{-\alpha} + \max \{t^\alpha_n, t_{n}^{2-\alpha}\}\left[   n^{-r(2\alpha+1)}
	\sum _{k=1} ^{\lceil n/2 \rceil -1}  k^{r(\alpha +1) -3} +  n^{-(2-\alpha)}\right]\\
&= C\Bigl(\frac{t_1}{t_n}\Bigr)^{2\alpha} t_n ^{\alpha} + \max \{t^\alpha_n, t_{n}^{2-\alpha}\}
	\left[ n^{-r(2\alpha+1)} \sum _{k=1} ^{\lceil n/2 \rceil -1}  k^{r(\alpha +1) -3} +  n^{-(2-\alpha)}\right]\\
&\leq C \max \{t^\alpha_n, t_{n}^{2-\alpha}\} \left[n^{-(2-\alpha)} + n^{-2r\alpha}\right],
\end{align*}
which is the desired result.
\end{proof}

\begin{lemma}\label{lem:CaputoTruncError}
There exists a constant $C$ such that
\begin{align*}
& \left| \left( L_{M,N} -( D_t^\alpha +{\cal L}) \right) u(x_m,t_n)\right|  \leq C
      \begin{cases}
     (t_n^\alpha n^{-\min \{2-\alpha, 2r\alpha\} } + h^2), & \text{ for } t_n\le 1 , \\
     (T^{2-\alpha} N^{-\min \{2-\alpha, 2r\alpha\} } + h^2), & \text{ for } t_n> 1 ,
    \end{cases}
\end{align*}

\end{lemma}
\begin{proof}
When $t_n \leq 1$ the result follows immediately from Lemmas~\ref{lem:decomposition}-\ref{lem:Remainder-CaputoTruncError}.
In the case where $t_n >1$, we simply note the following: when  $r> (2-\alpha)/(2\alpha)$, from Lemma~\ref{lem:Remainder-CaputoTruncError}, one has
\begin{eqnarray*}
\left| (D^\alpha_N -D^\alpha_t) v(x_m,t_n)\right|
     &\leq & Ct^{2-\alpha}_n  n^{-(2-\alpha)} \\
		&=& CT^{2-\alpha} N^{-(2-\alpha)}\Bigl(\frac{n}{N}\Bigr)^{(r-1)(2-\alpha)} \\
		&\leq& CT^{2-\alpha} N^{-(2-\alpha)}
\end{eqnarray*}
 and,  when  $r \leq (2-\alpha)/(2\alpha)$, one has
\begin{align*}
\left| (D^\alpha_N-D^\alpha_t) v(x_m,t_n)\right|
     &\le C t^{2-\alpha}_n   n^{-2r\alpha} \\
		&=  CT^{2-\alpha} N^{-2r\alpha} \left (\frac{n}{N} \right)^{(2-\alpha)r -2r\alpha}\le C T^{2-\alpha} N^{-2r\alpha}.
\end{align*}
\end{proof}

\subsection{Stability and error estimates}\label{sec:staberror}

This section gives error estimates for the finite difference scheme~\eqref{scheme}. They are obtained from the truncation error estimate of Lemma~\ref{lem:CaputoTruncError} and a discrete comparison principle together with some suitable barrier functions.

\begin{lemma}\label{lem:maxprin}
The scheme~\eqref{scheme} satisfies a discrete maximum principle, viz., if $f\ge 0$ and $\phi \ge 0$ in~\eqref{scheme}, then the discrete solution satisfies $u_m^n \ge 0$ for all~$m$ and~$n$.
\end{lemma}
\begin{proof}
In the scheme~\eqref{scheme}, the formulation~\eqref{CaputoN1} of the approximation to the Caputo time fractional derivative can be written as
$$
D^\alpha_N u^n_m =\sum^n_{k=0} \Theta_{nk} u^k_m,
$$
where
\begin{align*}
\Theta_{nk}
    & := -\frac{\alpha}{\Gamma(1-\alpha)} \int^{t_{k+1}}_{s=t_{k-1}} (t_n-s)^{-\alpha-1} \phi_k(s) \, ds, \quad 0 <k<n, \\
\Theta_{n0}
    & := -\frac1{t_n^\alpha \Gamma(1-\alpha)}-\frac{\alpha}{\Gamma(1-\alpha)} \int^{t_1}_{s=0} (t_n-s)^{-\alpha-1} \phi_0(s) \, ds, \\
\Theta_{nn}
    &  := - \sum^{n-1}_{k=0} \Theta_{nk},
\end{align*}
and $\phi_k$ are the basis functions~\eqref{basis}. The weights $\Theta_{nk}$ satisfy
$$
\Theta_{nn}>0 \ \text{ and } \Theta_{nk}<0 \ \text{ for } 0\le k <n,
$$
because of the positivity of the basis functions~\eqref{basis}.  Recalling also the discretisation of the spatial differential operator ${\cal L}$, we see that the matrix associated with the discrete operator $L_{M,N}$ is irreducibly diagonally dominant with positive diagonal entries  and non-positive  offdiagonal entries.
Hence the system matrix is an $M$-matrix and the result follows.
\end{proof}

The discrete maximum principle of Lemma 4 implies a discrete comparison principle, which can be exploited using discrete barrier functions, as is well known.
The next lemma gives a very useful criterion for constructing discrete barrier functions.

\begin{lemma}\label{lem:BF1}
If a grid function $b(x_m,t_n)$ satisfies
\begin{equation}\label{CondBarrierFunction}
b(x_m,0)=0 \quad \text{ and } \quad b(x_m,t_i) \le b(x_m,t_j) \text{ if } i\le j,
\end{equation}
then
\begin{equation}\label{LowBoundCaputo}
D_N^\alpha b(x_m,t_n) \ge \frac1{\Gamma(1-\alpha)}\frac{b(x_m,t_n)}{t^\alpha_n}\,.
\end{equation}
\end{lemma}
\begin{proof}
This result follows from the positivity of the basis functions $\phi_k$ in the derivation of~\eqref{CaputoN1}.
\end{proof}

The next theorem, which is the main result of the paper, gives error estimates for the fitted scheme~\eqref{scheme}.

\begin{theorem} \label{th:cgce}
The solution $u_m^n$ of the scheme~\eqref{scheme} satisfies
$$
|u(x_m,t_n)-u_m^n|
     \le C \max\left\{T^{2\alpha},T^{2-\alpha}\right\}  N^{-\min\{ 2-\alpha, 2r\alpha \} }+C h^2.
$$
\end{theorem}

\begin{proof}
We consider two cases.

\emph{Case (i):} Assume $2 r \alpha > 2-\alpha$. From Lemma~\ref{lem:CaputoTruncError}, one has the following estimate for the truncation error
\begin{align*}
&\hspace{-10mm}\left| \bigl(L_{M,N} -( D_t^\alpha +{\cal L}) \bigr) u(x_m,t_n)\right|  \\
	&\le  C t^{\alpha-\frac{2-\alpha}{r}}_n  \left[t^{\frac{2-\alpha}{r}}_n n^{-(2-\alpha)}\right]
		+ CT^{2-\alpha} N^{-(2-\alpha)} +C h^2  \\
	&=  C \left[t^{\alpha-\frac{2-\alpha}{r}}_n T^{(2-\alpha)/r} + T^{2-\alpha} \right] N^{-(2-\alpha)} +C h^2.
\end{align*}
Consider the barrier function
$$
C \left[b_1(x_m,t_n)+ \left(T^{2-\alpha}  N^{-(2-\alpha)}+ h^2\right)\Phi(x_m,t_n) \right],
$$
where $C$ is a sufficiently large positive constant. The grid function $b_1\ge 0$ is defined by
$$
b_1(x_m,t_n):= t^{2\alpha-\frac{2-\alpha}{r}}_n  T^{\frac{2-\alpha}{r}} N^{-(2-\alpha)},
$$
and the grid function $\Phi$ is given by
\begin{equation} \label{SpatialBarrierFunction}
\Phi(x_m,t_n):=  \frac{x_m (l-x_m)}{p\, l^2}.
\end{equation}
Observe that $b_1$ satisfies~\eqref{CondBarrierFunction} and the function $\Phi$ satisfies
\[
0\le  \Phi(x_m,t_n)  \le  \frac1{4p}\,, \quad L_{M,N} \Phi(x_m,t_n) \ge 2p.
\]
 Then Lemmas~\ref{lem:maxprin} and \ref{lem:BF1} give
 \begin{align*}
|u(x_m,t_n)-u_m^n|
    & \le C\left[b_1(x_m,t_n) + (T^{2-\alpha}  N^{-(2-\alpha)}+ h^2) \Phi(x_m,t_n)\right] \\
    & \le C t^{2\alpha-\frac{2-\alpha}{r}}_n T^{\frac{2-\alpha}{r}}  N^{-(2-\alpha)} + C (T^{2-\alpha}  N^{-(2-\alpha)}+ h^2) \\
    & = C t_n^{2\alpha}\left(\frac{T}{t_n}\right)^{\frac{2-\alpha}{r}} N^{-(2-\alpha)} + C (T^{2-\alpha}  N^{-(2-\alpha)}+ h^2)  \\
    & \le C  (T^{2-\alpha} +T^{2\alpha}) N^{-(2-\alpha)}+C h^2.
\end{align*}

\emph{Case (ii):} Assume $2 r \alpha \le 2-\alpha$. Then Lemma~\ref{lem:CaputoTruncError} gives
\begin{align*}
\left| \bigl(L_{M,N} -( D_t^\alpha +{\cal L}) \bigr) u(x_m,t_n)\right|
	&\le C (t^{\alpha}_n n^{-2r\alpha}   + T^{2-\alpha}  N^{-2r\alpha} +h^2) \\
	&\leq C \left[(T^{2\alpha} t^{-\alpha}_n   + T^{2-\alpha})  N^{-2r\alpha} +h^2\right]. 
\end{align*}
We now use the barrier function
$$
C \left[b_2(x_n,t_m)+(T^{2-\alpha}  N^{-2r\alpha} +h^2) \Phi(x_m,t_n) \right],
$$
where $C$ is a sufficiently large positive constant, $\Phi$ is given in~\eqref{SpatialBarrierFunction} and  $b_2$ is defined by
$$
b_2(x_m,t_n):= T^{2\alpha} N^{-2r\alpha}(1-e^{-t_n/t_1}).
$$
Observe that $b_2$ also satisfies~\eqref{CondBarrierFunction}. Furthermore, from Lemma~\ref{lem:BF1} we see that
\begin{align*}
D^\alpha_N b_2(x_m,t_n)
     \ge \frac1{\Gamma(1-\alpha)} \frac{b_2(x_m,t_n)}{t^\alpha_n}
     \ge C T^{2\alpha} \frac{1-e^{-1}}{\Gamma(1-\alpha)} t^{-\alpha}_n   N^{-2r\alpha}
\end{align*}
for $n\ge 1$.
Lemmas~\ref{lem:maxprin} and \ref{lem:BF1}   now yield
\begin{align*}
|u(x_m,t_n)-u_m^n|
    & \le C \left[b_2(x_m,t_n)+(T^{2-\alpha}  N^{-2r\alpha} +h^2) \Phi(x_m,t_n) \right] \\
    & \le C (T^{2\alpha}+T^{2-\alpha})  N^{-2r\alpha}+ C h^2.
\end{align*}
\end{proof}

In the next two remarks we discuss the error estimates given in Theorem~\ref{th:cgce}.

\begin{remark} [Error estimates for the optimal choice of the grading exponent $r$] \label{rem:ErrorEstimates}
 Theorem~\ref{th:cgce} shows that the scheme~\eqref{scheme} converges with order $2r\alpha$ if the grading exponent is chosen such that $r\le (2-\alpha)/(2\alpha)$ and, otherwise, it converges with order $2-\alpha$.

We propose to take  $r=\max\{1, (2-\alpha)/(2\alpha) \}$ as the optimal choice of the grading exponent; then the scheme~\eqref{scheme} converges with order $2-\alpha$ for all $0<\alpha<1$ and the clustering of too many grid points near $t=0$ is avoided.
Thus if $\alpha \ge 2/3$,  then  the mesh is uniform and Theorem~\ref{th:cgce} gives
\begin{align*}
|u(x_m,t_n)-u_m^n|       \le C T^{2\alpha}  N^{-(2-\alpha)}+C h^2. 
\end{align*}
On the other hand, if $\alpha < 2/3$, then  the mesh is graded and 
\begin{align*}
|u(x_m,t_n)-u_m^n|
  \le C T^{2-\alpha}  N^{-(2-\alpha)}+C h^2. 
\end{align*}

These error estimates, besides giving the order of convergence of the finite difference scheme~\eqref{scheme}, show the dependence on the final time $T$.  We shall illustrate in our numerical experiments that these estimates are sharp (see Tables~\ref{tb:Fittedgraded}, \ref{tb:FittedgradedT10} and~\ref{tb:Rates}).
\end{remark}

\begin{remark} [Error estimates for a suboptimal choice of $r$] \label{rem:ErrorEstimatesSpecial}
The error estimates for $\alpha < 2/3$ and $1\le r < (2-\alpha)/(2\alpha)$---this includes the approximation on a uniform mesh---are not discussed in the previous remark.

In this case, Theorem~\ref{th:cgce} shows that the time error discretization is bounded by $CT^{2-\alpha}N^{-2r\alpha}$. We have observed numerically that the order of convergence is $N^{-2r\alpha}$, as expected. But for some values of $r$ (e.g., $r=1$), we have also observed that the dependence on $T$ is $O(T^{2\alpha})$, not $T^{2-\alpha}$.

Some new error estimates are given in this remark that explain this behaviour. If one considers the barrier function
\begin{align*}
  C \left(b_2(x_n,t_m)+h^2 \Phi(x_m,t_n) + T^{2\alpha} N^{-2r\alpha} t_n^{2(1-\alpha)}\right)
\end{align*}
one obtains the error estimates
\begin{align*}
|u(x_m,t_n)-u_m^n|
    & \le C \left(T^{2\alpha} N^{-2r\alpha}+ h^2+T^{2\alpha} N^{-2r\alpha} t_n^{2(1-\alpha)} \right) \\
    & \le C T^{2\alpha} N^{-2r\alpha} \left(1+ t_n^{2(1-\alpha)}\right) + C h^2,
\end{align*}
and therefore
\begin{align*}
|u(x_m,t_n)-u_m^n| \le  C T^{2\alpha} N^{-2r\alpha} + C h^2, \quad \text {for } t_n \le 1. 
\end{align*}
Thus, if the maximum error occurs near  $t=0$ the growth rate is $T^{2\alpha}$ but not $T^{2-\alpha}$. This new error estimate is
illustrated in the numerical section in Tables~\ref{tb:Fitteduniform}, \ref{tb:FitteduniformT10} and~\ref{tb:Ratesuniform}.
\end{remark}
\section{Numerical results} \label{sec:Numerical}
In this section we  first give numerical results for two examples with known and unknown exact solutions that exhibit typical behaviour, i.e., their derivatives behave like~\eqref{Regu}. The examples are approximated with a standard scheme (which uses the classical L1 approximation for the Caputo time-fractional derivative and central difference approximation for the spatial variable; see for example~\cite{SORG17}) and the fitted scheme~\eqref{scheme}. Our results show that the fitted scheme produces  more accurate approximations  to the solution of both examples. The theoretical orders of convergence (TOC) from Theorem~\ref{th:cgce} and~\cite[Theorem 5.3]{SORG17} are included in the tables for comparison purposes.

The section finishes with an example whose exact solution is a smooth function. For this atypical example we find that on a uniform mesh our fitted method converges as the same rate as the classical L1 approximation.

\subsection{Numerical examples} 

\begin{example}\label{example1}
Consider the following problem:
$$
D_t^\alpha u- \, \frac{\partial^2 u}{\partial x^2}=f(x,t) \quad\text{for } (x,t)\in(0,\pi)\times(0,1],
$$
with initial condition $u(x,0)=\sin x$ for $0 < x < \pi$ and boundary conditions $u(0,t) = u(\pi,t)=0$ for $0\le t\le 1$.
The function $f$ is chosen such that the exact solution of the problem is
$$
u(x,t)=0.5 \left[ E_{\alpha,1}(-t^\alpha)+\cos(\pi t /3) \right] \sin x,
$$
which has a typical weak singularity at $t=0$  (see~\cite{SORG17}).
\end{example}

Since the exact solution of Example~\ref{example1} is known,   we can compute the exact errors in all the numerical approximations. The maximum nodal errors are denoted by
$$
E^{M,N}:=  \max_{(x_m,t_n)\in \bar Q}  \vert u(x_m,t_n)-u^n_m \vert
$$
and the corresponding  rate of convergence are computed using the usual formula
$$
\text{rate}_{M,N}:=  \log_2\left (\frac{E^{M,N}}{E^{2M,2N}} \right).
$$

 In~\cite[Theorem 5.3]{SORG17} it is established  that the standard L1 scheme on a graded mesh satisfies the error estimate
\begin{equation} \label{errorL1}
\max_{(x_m,t_n)\in \bar Q}  |u(x_m,t_n)-u_m^n|  \le C  \left( h^2 + N^{-\min\{2-\alpha, \,r\alpha  \}}  \right)
\end{equation}
for some constant $C$. Thus if $r=1$ (uniform mesh) the L1 scheme converges in time with order $\alpha$, and if $r\ge (2-\alpha)/\alpha$ it converges with order $2-\alpha$. The numerical results computed with this scheme for $r=1$ and $r=(2-\alpha)/\alpha$ are given in Tables~\ref{tb:L1uniformExact} and~\ref{tb:L1gradedExact} and they are in agreement with~\eqref{errorL1}.  In the case of a uniform mesh we note the slow rates of convergence for  small values of $\alpha$.

\begin{table}[h]
\caption{Standard L1 scheme: Maximum nodal errors and orders of convergence for Example~\ref{example1} with a uniform mesh and final time $T=1$; TOC come from~\eqref{errorL1}}
\begin{center}{ \label{tb:L1uniformExact}
\begin{tabular}{|c|c||c|c|c|c|c|}
 \hline  &TOC&  N=M=64 &N=M=128 & N=M=256 & N=M=512 & N=M=1024 \\
\hline $\alpha=0.2$ & 0.2
&1.855E-2 &1.729E-2 &1.603E-2 &1.477E-2 &1.353E-2 \\
&&0.101&0.109&0.118&0.127&
\\ \hline $\alpha=0.4$ & 0.4
&1.781E-2 &1.420E-2 &1.120E-2 &8.747E-3 &6.784E-3 \\
&&0.327&0.343&0.356&0.367&
\\ \hline $\alpha=0.6$  & 0.6
&9.117E-3 &6.119E-3 &4.081E-3 &2.711E-3 &1.797E-3 \\
&&0.575&0.584&0.590&0.594&
\\ \hline $\alpha=0.8$& 0.8
&3.067E-3 &1.781E-3 &1.026E-3 &5.901E-4 &3.391E-4 \\
&&0.785&0.795&0.798&0.799&
\\ \hline \hline
\end{tabular}}
\end{center}
\end{table}

\begin{table}[h]
\caption{Standard L1 scheme: Maximum nodal errors and orders of convergence for Example~\ref{example1} for a graded mesh with an optimal  grading of $r=(2-\alpha)/\alpha$ and final time $T=1$; TOC come from~\eqref{errorL1}}
\begin{center}{ \label{tb:L1gradedExact}
\begin{tabular}{|c|c||c|c|c|c|c|}
 \hline &TOC & N=M=64 &N=M=128 & N=M=256 & N=M=512 & N=M=1024 \\
\hline $\alpha=0.2$ & 1.8
&1.933E-4 &6.570E-5 &2.137E-5 &6.761E-6 &2.097E-6 \\
&&1.557&1.620&1.661&1.689&
\\ \hline $\alpha=0.4$ & 1.6
&3.673E-4 &1.262E-4 &4.308E-5 &1.463E-5 &4.948E-6 \\
&&1.541&1.551&1.558&1.564&
\\ \hline $\alpha=0.6$ & 1.4
&7.097E-4 &2.785E-4 &1.083E-4 &4.181E-5 &1.605E-5 \\
&&1.349&1.363&1.373&1.381&
\\ \hline $\alpha=0.8$ & 1.2
&1.099E-3 &5.127E-4 &2.366E-4 &1.081E-4 &4.899E-5 \\
&&1.100&1.115&1.130&1.142&
\\ \hline \hline
\end{tabular}}
\end{center}
\end{table}

We next use  the fitted scheme~\eqref{scheme} to generate approximate solutions of  Example~\ref{example1}. In Tables~\ref{tb:Fitteduniform} and~\ref{tb:Fittedgraded} the maximum nodal errors and the orders of convergence are given with $r=1$  and $r=\max\{1, (2-\alpha)/(2\alpha) \}$, respectively. Note that when $r=1$, Theorem~\ref{th:cgce} indicates that the order is $\min \{2\alpha, 2-\alpha \}$ and by Remark~\ref{rem:ErrorEstimates} the order is $2-\alpha$ if the optimal grading is used. The computed orders of convergence are in agreement with Theorem~\ref{th:cgce} and Remark~\ref{rem:ErrorEstimates}. Comparing the numerical results for both schemes, we observe that the fitted scheme is more accurate  both on uniform and graded meshes.

\begin{table}[h]
\caption{Fitted scheme: Maximum nodal errors and orders of convergence for Example~\ref{example1} on a uniform mesh and final time $T=1$; TOC come from Theorem~\ref{th:cgce}}
\begin{center}{ \label{tb:Fitteduniform}
\begin{tabular}{|c|c||c|c|c|c|c|}
 \hline  & TOC & N=M=64 &N=M=128 & N=M=256 & N=M=512 & N=M=1024 \\
 \hline $\alpha=0.2$ & 0.4
&1.971E-3 &1.656E-3 &1.392E-3 &1.161E-3 &9.591E-4 \\
&&0.251&0.251&0.262&0.275&
\\ \hline $\alpha=0.4$ & 0.8
&2.207E-3 &1.410E-3 &8.808E-4 &5.401E-4 &3.263E-4 \\
&&0.647&0.678&0.705&0.727&
\\ \hline $\alpha=0.6$ & 1.2
&9.401E-4 &4.409E-4 &2.019E-4 &9.090E-5 &4.081E-5 \\
&&1.092&1.127&1.151&1.155&
\\ \hline $\alpha=0.8$ & 1.2
&3.435E-4 &1.515E-4 &6.670E-5 &2.927E-5 &1.281E-5 \\
&&1.181&1.184&1.188&1.192&
\\ \hline \hline
\end{tabular}}
\end{center}
\end{table}

\begin{table}[h]
\caption{Fitted scheme: Maximum nodal errors and orders of convergence for Example~\ref{example1} on a graded mesh  with an optimal grading of $r=\max\{1, (2-\alpha)/(2\alpha)\}$ and final time $T=1$; TOC come from Theorem~\ref{th:cgce}}
\begin{center}{ \label{tb:Fittedgraded}
\begin{tabular}{|c|c||c|c|c|c|c|}
 \hline & TOC & N=M=64 &N=M=128 & N=M=256 & N=M=512 & N=M=1024 \\
\hline $\alpha=0.2$ & 1.8
&1.451E-4 &4.891E-5 &1.584E-5 &4.996E-6 &1.546E-6 \\
&&1.568&1.626&1.665&1.692&
\\ \hline $\alpha=0.4$ & 1.6
&1.366E-4 &4.639E-5 &1.657E-5 &5.883E-6 &2.045E-6 \\
&&1.558&1.486&1.494&1.524&
\\ \hline $\alpha=0.6$ & 1.4
&4.762E-4 &1.965E-4 &7.986E-5 &3.202E-5 &1.270E-5 \\
&&1.277&1.299&1.319&1.335&
\\ \hline $\alpha=0.8$ & 1.2
&3.435E-4 &1.515E-4 &6.670E-5 &2.927E-5 &1.281E-5 \\
&&1.181&1.184&1.188&1.192&
\\ \hline \hline
\end{tabular}}
\end{center}
\end{table}

We now illustrate the dependence of the errors on the final time $T$.  Consider again Example~\ref{example1} but with final time $T=10$ instead of $T=1$. The maximum nodal errors for the fitted scheme using  a graded mesh with $r=\max\{1, (2-\alpha)/(2\alpha) \}$ are displayed in Table~\ref{tb:FittedgradedT10}; they are clearly greater than in Table~\ref{tb:Fittedgraded}, so there is some dependence on the final time~$T$. In Table~\ref{tb:Rates} we show the computed ``growth rates"
$$
\log_{10}(E^{M,N}_{T=10}/E^{M,N}_{T=1}),
$$
where $E^{M,N}_{T=1}$ and $E^{M,N}_{T=10}$ are the maximum errors for $T=1$ and $T=10$ (they are given in Tables~\ref{tb:Fittedgraded} and~\ref{tb:FittedgradedT10}, respectively). These numerical results are in agreement with
Remark~\ref{rem:ErrorEstimates} showing that  the  dependence on $T$ is approximately $O(T^{2-\alpha})$ if $\alpha<2/3$ and $O(T^{2\alpha})$ if $\alpha \ge 2/3$.

\begin{table}[h]
\caption{Fitted scheme: Maximum nodal errors and orders of convergence for Example~\ref{example1} on a graded mesh  with an optimal grading of $r=\max\{1, (2-\alpha)/(2\alpha) \}$ and final time $T=10$; TOC come from Theorem~\ref{th:cgce}}
\begin{center}{ \label{tb:FittedgradedT10}
\begin{tabular}{|c|c||c|c|c|c|c|}
 \hline &  TOC & N=M=64 &N=M=128 & N=M=256 & N=M=512 & N=M=1024 \\
 \hline $\alpha=0.2$ & 1.8
&8.058E-3 &2.699E-3 &8.511E-4 &2.655E-4 &8.145E-5 \\
&&1.578&1.665&1.681&1.705&
\\ \hline $\alpha=0.4$ & 1.6
&6.509E-3 &2.285E-3 &7.877E-4 &2.683E-4 &9.064E-5 \\
&&1.510&1.536&1.554&1.566&
\\ \hline $\alpha=0.6$ & 1.4
&7.169E-3 &2.797E-3 &1.080E-3 &4.138E-4 &1.687E-4 \\
&&1.358&1.373&1.384&1.295&
\\ \hline $\alpha=0.8$ & 1.2
&1.434E-2 &6.374E-3 &2.802E-3 &1.226E-3 &5.352E-4 \\
&&1.170&1.186&1.193&1.196&
\\ \hline \hline
\end{tabular}}
\end{center}
\end{table}

\begin{table}[h]
\caption{Fitted scheme: Growth rates for Example~\ref{example1} on a graded mesh  with an optimal grading of $r=\max\{1, (2-\alpha)/(2\alpha) \}$ and final times $T=1,10$}
\begin{center}{ \label{tb:Rates}
\begin{tabular}{|c|c||c|c|c|c|c|}
 \hline &   N=M=64 &N=M=128 & N=M=256 & N=M=512 & N=M=1024 \\
 \hline $\alpha=0.2$
 &  1.745 &   1.742  &  1.730 &   1.725 &   1.722
 \\
 \hline $\alpha=0.4$
 &   1.678  &  1.693  &  1.677  &  1.659 &   1.647
 \\
 \hline $\alpha=0.6$
 &  1.178 &   1.153 &   1.131 &   1.111 &   1.123
 \\
 \hline $\alpha=0.8$
 &  1.621 &   1.624  &  1.623  &  1.622  &  1.621
 \\
 \hline \hline
\end{tabular}}
\end{center}
\end{table}

The numerical results for the fitted scheme using  a uniform mesh for $T=10$ are displayed in Table~\ref{tb:FitteduniformT10}, which are again greater than for $T=1$ (see Table~\ref{tb:Fitteduniform}). In Table~\ref{tb:Ratesuniform} we show the computed growth rates from the maximum errors for $T=1$ (Table~\ref{tb:Fitteduniform}) and $T=10$ (Table~\ref{tb:FitteduniformT10}). These numerical results are in agrement with
Remark~\ref{rem:ErrorEstimates} for $\alpha \ge 2/3$ and Remark~\ref{rem:ErrorEstimatesSpecial} for $\alpha < 2/3$, showing that  the dependence on $T$ is approximately $O(T^{2\alpha})$ for all the values of $\alpha$ considered in these tables.

\begin{table}[h]
\caption{Fitted scheme: Maximum nodal errors and orders of convergence for Example~\ref{example1} on a uniform mesh and final time $T=10$; TOC come from Theorem~\ref{th:cgce}}
\begin{center}{ \label{tb:FitteduniformT10}
\begin{tabular}{|c|c||c|c|c|c|c|}
 \hline & TOC & N=M=64 &N=M=128 & N=M=256 & N=M=512 & N=M=1024 \\
 \hline $\alpha=0.2$ & 0.4
&1.908E-3 &2.370E-3 &2.259E-3 &2.006E-3 &1.731E-3 \\
&&-0.313&0.069&0.172&0.213&
\\ \hline $\alpha=0.4$ & 0.8
&5.182E-3 &4.834E-3 &3.624E-3 &2.485E-3 &1.622E-3 \\
&&0.100&0.416&0.544&0.616&
\\ \hline $\alpha=0.6$ & 1.2
&6.184E-3 &3.987E-3 &2.276E-3 &1.164E-3 &5.594E-4 \\
&&0.633&0.808&0.968&1.056&
\\ \hline $\alpha=0.8$ & 1.2
&1.434E-2 &6.374E-3 &2.802E-3 &1.226E-3 &5.352E-4 \\
&&1.170&1.186&1.193&1.196&
\\ \hline \hline
\end{tabular}}
\end{center}
\end{table}

\begin{table}[h]
\caption{Fitted scheme: Growth rates for Example~\ref{example1} on a uniform mesh  and final times $T=1,10$}
\begin{center}{\label{tb:Ratesuniform}
\begin{tabular}{|c||c|c|c|c|c|}
 \hline  &  N=M=64 &N=M=128 & N=M=256 & N=M=512 & N=M=1024 \\
  \hline $\alpha=0.2$
  &-0.014    &0.156    &0.210    &0.238    &0.256
 \\
  \hline $\alpha=0.4$
 &0.371    &0.535    &0.614    &0.663    &0.696
    \\
  \hline $\alpha=0.6$
         &0.818    &0.956    &1.052    &1.107    &1.137
 \\
 \hline $\alpha=0.8$
   &  1.621 &   1.624  &  1.623  &  1.622  &  1.621
 \\ \hline \hline
\end{tabular}}
\end{center}
\end{table}

\begin{example}\label{example2}

Consider the  test problem
\begin{align*}
D_t^\alpha u- \frac{\partial^2 u}{\partial x^2} + (1+x) u & = x(\pi-x) (1+t^4)+t^2 \\
\intertext{for $(x,t)\in Q:=(0,\pi)\times(0,1]$, with}
u(0,t)=0, \quad u(\pi,t) =0, & \quad \text{for }t\in(0,1],  \\
u(x,0)=\sin x, & \quad \text{for } x\in[0,\pi].
\end{align*}
\end{example}

The exact solution of Example~\ref{example2} is unknown and the orders of convergence in the computed solutions are estimated using the two-mesh principle~\cite{FHMORS00}. The maximum nodal two-mesh differences are defined by
\[
D^{M,N}:=\max_{0\le m\le M} \max_{0\le n\le N} \vert u^n_{m}-z^{2n}_{2m}\vert,
\]
where $z^n_m$, with $0\le m \le 2M$ and $0\le n \le 2N$ is the computed solution by the scheme~\eqref{scheme} which uses a uniform in the spatial direction with $2M$ mesh intervals, and a graded mesh
\[
t_n= \big(n/(2N)\big)^r \quad\text{for } 0\le n \le 2N,
\]
in the temporal direction. The maximum nodal two-mesh differences are used to compute the estimated rate of convergence in a standard way
\[
\log_2\left (\frac{D^{M,N}}{D^{2M,2N}} \right).
\]
The maximum nodal two-mesh differences and their corresponding estimated orders of convergence for Example~\ref{example2} with the fitted scheme are given in Table~\ref{tb:FittedGradedUnknown}, where the grading exponent is $r=\max\{1, (2-\alpha)/(2\alpha) \}$ (see Remark~\ref{rem:ErrorEstimates}). The numerical results are in agreement with the theoretical rates of convergence established in Theorem~\ref{th:cgce}.

\begin{table}[h]
\caption{Fitted scheme: Maximum nodal two-mesh differences and orders of convergence for Example~\ref{example2} on a graded mesh with an optimal  grading of  $r=\max\{1, (2-\alpha)/(2\alpha) \}$; TOC come from Theorem~\ref{th:cgce}}
\begin{center}{ \label{tb:FittedGradedUnknown}
\begin{tabular}{|c||c|c|c|c|c|c|}
 \hline  & TOC &  N=M=64 &N=M=128 & N=M=256 & N=M=512 & N=M=1024 \\
\hline $\alpha=0.2$ & 1.8
&7.411E-4 &2.438E-4 &7.749E-5 &2.408E-5 &7.365E-6 \\
&&1.604&1.653&1.686&1.709&
\\ \hline $\alpha=0.4$ & 1.6
&6.979E-4 &2.712E-4 &1.030E-4 &3.813E-5 &1.376E-5 \\
&&1.364&1.397&1.433&1.470&
\\ \hline $\alpha=0.6$ & 1.4
&2.105E-3 &9.965E-4 &4.606E-4 &2.006E-4 &8.506E-5 \\
&&1.079&1.113&1.200&1.237&
\\ \hline $\alpha=0.8$ & 1.2
&1.866E-3 &8.863E-4 &4.056E-4 &1.820E-4 &8.062E-5 \\
&&1.074&1.128&1.156&1.175&
\\ \hline \hline
\end{tabular}}
\end{center}
\end{table}

Finally, we consider an example where the solution is smooth (i.e., all derivatives needed in the numerical analysis are bounded on the closed domain), since this class of problem has been considered by several researchers (but see~\cite{StyFCAA}). For such problems,  our earlier analysis can be modified to show that on a uniform mesh the fitted scheme satisfies the error bound
\begin{equation}\label{cgcesmooth}
\vert u(x_m,t_n) - u^n_m \vert \leq CT^{2-\alpha} N ^{-(2-\alpha)} +Ch^2
\end{equation}
for all $m,n$.

\begin{example}\label{example3}
Take $p=1, c(x)=x^2$ and choose the solution of \eqref{prob} to be
$$
u(x,t)=(1+t^3)(4 x (1-x))^2+5 t^{3+\alpha}, \quad (x,t)\in [0,1]\times[0,1].
$$
Then the right-hand side $f$, the initial condition and the boundary conditions are determined by the statement of~\eqref{prob} and our choices of $p, c$ and $u$.
\end{example}

Example~\ref{example3} is solved using our fitted scheme on a uniform mesh. The maximum errors and orders of convergence are given in Table~\ref{tb:smooth}, where we have taken $M=4N$  in order that the error from the temporal discretisation of the fractional term dominates the spatial error. It is clear that the method converges with order $2-\alpha$, in agreement with~\eqref{cgcesmooth}.

\begin{table}[h]
\caption{Fitted scheme on a uniform mesh: Maximum nodal errors and orders of convergence for an example with smooth solution; TOC come from~\eqref{cgcesmooth}}
\begin{center}{ \label{tb:smooth}
\begin{tabular}{|c||c|c|c|c|c|c|}
 \hline  & TOC &  M=64 &M=128 & M=256 & M=512 & M=1024 \\
  &  &  N=16 &N=32 & N=64 & N=128 & N=256 \\
\hline $\alpha=0.2$ & 1.8
& 4.832E-3 &1.407E-3 &4.086E-4 &1.184E-4 &3.425E-5\\
&& 1.780&1.784&1.787&1.789&
\\ \hline $\alpha=0.4$ & 1.6
& 1.323E-2 &4.451E-3 &1.488E-3 &4.955E-4 &1.645E-4
\\
&& 1.571&1.581&1.587&1.590&
\\ \hline $\alpha=0.6$ & 1.4
& 3.429E-2 &1.326E-2 &5.083E-3 &1.940E-3 &7.385E-4
\\
&& 1.371&1.383&1.390&1.394&
\\ \hline $\alpha=0.8$ & 1.2
& 8.410E-2 &3.726E-2 &1.638E-2 &7.167E-3 &3.129E-3
\\
&&1.174&1.186&1.192&1.196&
\\ \hline \hline
\end{tabular}}
\end{center}
\end{table}

\bibliographystyle{spmpsci}      

%
%
\section*{Acknowledgements}
The research of Jos\'e Luis Gracia was partly supported by the Institute of Mathematics and Applications (IUMA), the project MTM2016-75139-R  and the Diputaci\'{o}n General de Arag\'{o}n. The research of Martin Stynes was supported in part by the National Natural Science Foundation of China under grants 91430216 and NSAF-U1530401.


%
%

\end{document}